\documentclass[oneside,a4paper,english, 11pt]{amsart}

\input{artmymacrosSurvey.sty}

\title{Serre weights, Galois deformation rings, and local models}

\author{Daniel Le}
\address{Department of Mathematics,
Purdue University,
150 N. University Street, 
West Lafayette, IN 47907-2067}
\email{ledt@purdue.edu}

\author{Bao V.~Le Hung}
\address{Department of Mathematics,
Northwestern University, 
2033 Sheridan Road, 
Evanston, Illinois 60208, USA}
\email{lhvietbao@googlemail.com}

\begin{document}

\begin{abstract} 
We survey some recent progress on generalizations of conjectures of Serre concerning the cohomology of arithmetic groups, focusing primarily on the ``weight" aspect. 
This is intimately related to (generalizations of) a conjecture of Breuil and M\'ezard relating the geometry of potentially semistable deformation rings to modular representation theory. 
Recently, B.~Levin, S.~Morra, and the authors established these conjectures in tame generic contexts by constructing projective varieties (local models) in mixed characteristic whose singularities model, in generic cases, those of tamely potentially crystalline Galois deformation rings for unramified extensions of $\Qp$ with small regular Hodge--Tate weights.
\end{abstract}

\maketitle


\section{Introduction}

The study of congruences between automorphic forms has a long and rich tradition. 
A paradigm shift occurred when Deligne's construction of Galois representations attached to classical holomorphic Hecke eigenforms opened the door to the study of congruences of automorphic forms through congruences of Galois representations. 
In fact, conjectures of Fontaine--Mazur--Langlands and Serre suggest that these are really two sides of the same coin. 

Let $p$ be a prime. 
Recall that Serre's conjecture asserts that every continuous, odd, and irreducible Galois representation $\rhobar: G_{\Q} \ra \GL_2(\overline{\F}_p)$ of the absolute Galois group of $\Q$ is \emph{modular}, i.e.~arises from the reduction of a Galois representation attached to a classical modular form. 
It furthermore asserts a refinement which specifies the minimal weight and level at which one may find such a modular form in terms of local properties of $\rhobar$. 
Turning the perspective around, if one begins with a modular $\rhobar$, then the refinement predicts congruences between modular forms of many different levels and weights. 

While the level part generalizes quite naturally, the weight part is subtler, because it turns out to be inextricably linked to integral $p$-adic Hodge theory. The goal of this survey is to describe the circle of ideas surrounding recent developments on the weight part of Serre's conjecture. 

\subsection{Overview} 
In \S \ref{sec:coh}, we articulate the main questions of interest in the context of cohomological automorphic forms, especially motivating the representation theoretic perspective on congruences. 

In \S \ref{sec:reptheory} and \ref{sec:galrep}, we briefly discuss background on modular representation theory and Galois representations. This will be necessary to state the higher dimensional generalizations of the weight part of Serre's conjecture. 

In \S \ref{sec:conj}, we state generalizations of Serre's conjecture due to Ash, Gee, Herzig, and Savitt among many others which provide conjectural answers to these questions in many cases. In \S \ref{sec:result}, we narrow our focus to two cases, definite unitary groups in joint work with B.~Levin and S.~Morra and $\GL_n$ over CM fields, where we established some conjectures of \cite{GHS} when $\rhobar$ is tame and sufficiently generic at $p$. 
In the proofs, the Kisin--Taylor--Wiles method plays a vital role in reducing a global problem to a local one. 
As the internal details of the Kisin--Taylor--Wiles method are orthogonal to our goals, we have chosen to axiomatize its essential output and explain how it is used. 

Finally, \S \ref {sec:LM} summarizes the key results on local models for local deformation rings which are the essential ingredients to prove the results on the weight part of Serre's conjecture. 

\subsection{Acknowledgments}

This survey is based on two talks given by the authors at the International Colloquium on Arithmetic Geometry at TIFR in January 2020. 
We thank the organizers for the invitation to this wonderful mathematical and cultural event. 


D.L.~was supported by the National Science Foundation under agreement No.~DMS-1703182. B.LH.~acknowledges support from the National Science Foundation under grant Nos.~DMS-1128155, DMS-1802037 and the Alfred P.~Sloan Foundation. 

\section{Cohomology of arithmetic manifolds} \label{sec:coh}

Let $\bG$ be a connected reductive group over $\Q$.
Let $A_\infty^\circ$ be the connected component of the group of $\R$-points of a maximal $\Q$-split torus in the center of $\bG$, and let $K_\infty^\circ$ be a maximal connected compact subgroup of $\bG(\R)$.
For a compact open subgroup $K_f \subset \bG(\A)$, consider the adelic double quotient
\[
Y(K_f) \defeq \bG(\Q)\backslash \bG(\A)/A_\infty^\circ K_\infty^\circ K_f.
\]
In various places, we require technical properties of $K_f$ that can always be attained by passing to a finite index subgroup. 
Furthermore, these required properties are preserved under passing to a finite index subgroup (away from a finite set of places). 
We assume throughout that $K_f$ is sufficiently small (cf.~\cite[\S 9.1]{MLM}). 
In particular, $K_f$ is neat so that $Y(K_f)$ is naturally a real manifold. 
Moreover, $Y(K_f)$ can be rewritten as a finite disjoint union of quotients of the symmetric space $\bG(\R)/A_\infty^\circ K_\infty^\circ$ by subgroups of $\bG(\Q)$ which are discrete and of finite covolume in $\bG(\R)$.
Finally, $Y(K_f)$ is homotopy equivalent to its Borel--Serre compactification (cf.~\cite{Borel-Serre}), a compact real manifold with corners and, in particular, a finite CW complex. 

\subsection{Rational cohomology} \label{sec:ratcoeff}

Let $V$ be an algebraic representation of $\bG$ over $\Q$. 
Then let $\cV_{\Q}$ be the $\Q$-local system 
\[
\bG(\Q) \Big\backslash (\bG(\A) /A_\infty^\circ K_\infty^\circ K_f) \times V(\Q),
\]
where $\bG(\Q)$ acts diagonally.
Then the (finite-dimensional) sheaf cohomology groups $H^*(Y(K_f),\cV_{\Q})$ have a convolution action by the double coset (Hecke) algebra $\Q[K_f\backslash \bG(\bA_f) /K_f]$. 
(We will consider sequences of cohomology groups as objects in appropriate bounded derived categories. In most instances, little is lost if $H^*(Y(K_f),\cV_{\Q})$ is replaced by $\oplus_{i\in \Z} H^i(Y(K_f),\cV_{\Q})$.)
The significance of these cohomology groups stems from the fact that the Hecke module $H^*(Y(K_f),\cV_{\Q}) \otimes_{\Q} \C$ can be computed by automorphic forms by work of  Matsushima and Franke \cite{matsushima,franke}. This is essentially Hodge theory for the locally symmetric manifolds $Y(K_f)$.
Suppose that $A_\infty^\circ$ acts on $V(\C)$ through a character. 
If we let $\cA(K_f)$ denote the space of automorphic forms on $Y(K_f)$, then 
\begin{equation}\label{eqn:matsushima}
H^*(Y(K_f),\cV_{\Q}) \otimes_{\Q} \C \cong H^*(\mathfrak{g},\mathfrak{p},(\cA(K_f) \otimes_{\C}V(\C))^{A_\infty^\circ}), 
\end{equation}
where $\mathfrak{g}$ is the Lie algebra of the group of real points of the intersection of all kernels of rational characters of $\bG$. 
Let $\cV_{\C}$ be the local system 
\[
(\bG(\Q)\backslash \bG(\A)/K_f) \times V(\C) \Big/ A_\infty^\circ K_\infty^\circ,
\]
where the right action of $A_\infty^\circ K_\infty^\circ$ is diagonal with the inverse of the left action on $V(\C)$.
Then the bijection
\begin{align}
\bG(\A) \times V(\C) &\ra \bG(\A) \times V(\C)\\
(g,v) &\mapsto (g,g_\infty^{-1}v)
\end{align}
induces an isomorphism $\cV_{\Q} \otimes_{\Q} \C \risom \cV_{\C}$.
This is the first step in establishing \eqref{eqn:matsushima}.

We will assume that we can write $K_f$ as $K_\Sigma K^\Sigma$ where $\Sigma$ is a finite set of finite places and $K^\Sigma = \prod_{\ell\notin \Sigma} K_\ell$ where for all $\ell\notin \Sigma$, $K_\ell$ is a hyperspecial subgroup of $\bG(\Q_\ell)$ (in particular, we assume that $\bG$ is unramified at places $\ell\notin \Sigma)$. 
Then $T^\Sigma_{\Q} \defeq \Q[K^\Sigma \backslash \bG(\A^\Sigma)/K^\Sigma]$ is commutative (see \S \ref{sec:satake}). Since $H^*(Y(K_f),\cV_{\Q})$ are finite dimensional $\Q$-vector spaces, the eigenvalues of the $T^\Sigma_{\Q}$-action on the part of $\cA(K_f)$ which contributes to $H^*(Y(K_f),\cV_{\C})$ (the \emph{cohomological automorphic forms}) are algebraic numbers.


\subsection{Classical modular forms} \label{sec:classical}
If $\bG = \GL_2$, then $Y(K_f)$ is a modular curve and has the additional structure of a variety defined over $\Q$. 
Any irreducible algebraic representation of $\GL_2$ is of the form $V(a,b) \defeq \Sym^{a-b} (\Q^2) \otimes_\Q \det^b$.
Let $T_\ell \in \Q[\GL_2(\Z_\ell) \backslash \GL_2(\Q_\ell)/\GL_2(\Z_\ell)]$ be the double coset operator \[\GL_2(\Z_\ell) \begin{pmatrix}
\ell & 0 \\
0 & 1
\end{pmatrix}
\GL_2(\Z_\ell).\]
A well-known incarnation of \eqref{eqn:matsushima} is that there is a normalized Hecke eigenform 
\[
f(z) = \sum_{n=0}^\infty a_n q^n, \qquad \textrm{with } q = e^{2\pi i z}
\]
(i.e.~$a_1=1$) of weight $k \geq 2$ and level $K_f$ (or $K_f \cap \SL_2(\Z)$) if and only if there is a $T^\Sigma_{\Q}$-eigenvector in $H^1(Y(K_f),\cV(b+k-2,b))$ such that $T_\ell$ acts by $\ell^b a_\ell$ for all $\ell \notin \Sigma$.

It is well-known that the space of modular forms has a basis with integral $q$-expansions whose $\Z$-span is Hecke stable.
In particular, $(a_\ell)_\ell$ are not just algebraic numbers, but are in fact algebraic integers. This gives one way to make the notion of congruences between eigenforms precise: one asks for a congruence between the (integral) Fourier coefficients.

It turns out that there are a lot of congruences between $q$-expansions of integral Hecke eigenforms.
A basic example comes from the Eisenstein series
\[
G_k(z) \defeq -\frac{B_k}{2k} + \sum_{n\geq 1} \sigma_{k-1}(n)q^n, \qquad k \geq 4 \textrm{ even},
\]
where $B_k$ is the $k$-th Bernoulli coefficient.
Fixing a (rational) prime $p$, the mod $p$ $q$-series $G_k \pmod p$ depends only on $k \pmod {p-1}$.
Another well-known example is the congruence
\begin{equation}\label{eqn:cuspcong}
\Delta(z) \defeq q\prod_{m=1}^\infty (1-q^n)^{24} \equiv q\prod_{m=1}^\infty (1-q^n)^2 (1-q^{11n})^2 \pmod {11}
\end{equation}
between the unique normalized cuspforms of level $\Gamma(1)$ and weight $12$ and of level $\Gamma_0(11)$ and weight $2$, respectively.

The above notion of congruences between eigenforms is essentially equivalent to congruences between the system of Hecke eigenvalues on rational cohomology, and thus can also be detected by contemplating the action of (suitably integral) Hecke operators on cohomology with integral coefficients. It turns out that this shift of perspective from $q$-expansion to integral cohomology (initiated by Ash--Stevens) will give a systematic mechanism to explain congruences between automorphic forms via representation theory.



\subsection{Integral structure}\label{sec:integral}

Fix a prime $p$ and suppose that $K_f$ factors as the product $K_f^p K_p$.
We fix an algebraic closure $\ovl{\Q}_p$ of $\Q_p$ and let $E$ be a subfield of $\ovl{\Q}_p$ of finite degree over $\Q_p$.
By replacing $E$ if necessary, we will assume that $E$ is sufficiently large.
Let $\cO$ be the ring of integers of $E$ with uniformizer $\varpi$ and $\F$ be the residue field.
We define $\cV_{E}$ to be the nonarchimedean analogue
\[
(\bG(\Q)\backslash \bG(\A)/A_\infty^\circ K_\infty^\circ K_f^p) \times V(E)\Big/ K_p
\]
of $\cV_{\C}$ in \S \ref{sec:ratcoeff}, where $K_p$ acts diagonally (using the natural right action on $\bG(\Q)\backslash \bG(\A)/A_\infty^\circ K_\infty^\circ K_f^p$ and the inverse of the natural left action on $V(E)$).
Then as before, the map
\begin{align}
\bG(\A) \times V(E) &\ra \bG(\A) \times V(E)\\
(g,v) &\mapsto (g,g_p^{-1}v)
\end{align}
induces an isomorphism $\cV_{\Q} \otimes_{\Q} E \risom \cV_{E}$.
As $K_p$ is a compact group, there exists a $K_p$-stable $\cO$-lattice $W$ in $V(E)$.
If we let 
\begin{equation}\label{eqn:Zplocsys}
\cW\defeq (\bG(\Q)\backslash \bG(\A)/A_\infty^\circ K_\infty^\circ K_f^p \times W)/K_p, 
\end{equation}
then the map 
\begin{equation} \label{eqn:integral}
H^*(Y(K_f),\cW) \ra H^*(Y(K_f),\cV_{E}) \cong H^*(Y(K_f),\cV_{\Q}) \otimes_{\Q} E
\end{equation}
gives a natural integral structure on $H^*(Y(K_f),\cV_{\Q}) \otimes_{\Q} E$. 
In fact, the definition \eqref{eqn:Zplocsys} makes sense for any $\cO[\![K_p]\!]$-module and defines a functor from $\cO[\![K_p]\!]$-modules to local systems on $Y(K_f)$.
We caution that \eqref{eqn:integral} may not be injective in any given degree.
Indeed, $H^*(Y(K_f),\cW)$ may contain torsion, and in fact this torsion is expected to be abundant and to play an important role in connecting cohomological automorphic forms and Galois representations.

\subsection{Congruences between Hecke eigensystems}\label{sec:congruences}

Let $T^\Sigma_{\cO}$ be the Hecke algebra $\cO[K^\Sigma\backslash \bG(\A^\Sigma)/K^\Sigma]$, which acts naturally on $H^*(Y(K_f),\cW)$.
As $H^*(Y(K_f),\cW)$ is a finite $\cO$-module, there are only finitely many maximal ideals $\fm \subset T^\Sigma_{\cO}$ for which the localization $H^*(Y(K_f),\cW)_{\fm}$ is nonzero. 
These localized modules record congruences between systems of Hecke eigenvalues:
a Hecke eigenclass $H^*(Y(K_f),\cW)[\frac{1}{p}]$ survives in the localization $H^*(Y(K_f),\cW)_{\fm}[\frac{1}{p}]$ if and only if its (automatically integral) system of Hecke eigenvalues lifts the mod $p$ system given by $\fm$.

However, for a fixed $\fm$, there may be various local systems $\cW$ for which $H^*(Y(K_f),\cW)_{\fm}$ is nonzero.
Indeed, we saw in \S \ref{sec:classical} that if $\fm_{G_k}$ corresponds to (the system of Hecke eigenvalues of) the Eisenstein series $G_k \pmod p$ for $4\leq k \leq p+1$ and $\cW(a,b)$ corresponds to the lattice $W(a,b) \defeq \Sym^{a-b} \Z_p^2 \otimes \det^b$ for $a \geq b$, then $H^*(Y(\GL_2(\widehat{\Z})),\cW(k'-2,0))_{\fm_{G_k}}$ is nonzero for all $k' \equiv k \pmod {p-1}$.
(Since $\GL_2(\widehat{\Z})$ is not neat, these cohomology groups should be interpreted as the cohomology groups of an orbifold.)
Furthermore, if $\fm_\Delta$ corresponds to the Ramanujan Delta function mod $11$, then both $H^*(Y(K_f),\Z_p)_{\fm_\Delta}$ and $H^*(Y(\GL_2(\widehat{\Z})),\cW(10,0))_{\fm_\Delta}$ are nonzero where $K_f$ corresponds to the congruence subgroup $\Gamma_0(11)$.

The upshot of our discussion above is that congruences between eigenforms can be thought as the non-vanishing of localized cohomology for many different coefficient sheaves.
Thus a complete classification of such congruences is equivalent to the following question:

\begin{ques}\label{ques:Serre1}
Given a mod $p$ Hecke eigensystem $\fm$, for which $\cO$-local systems $\cW$ on $Y(K_f)$ is $H^*(Y(K_f),\cW)_{\fm}$ nonzero?
\end{ques}

Serre studied this question extensively in the case of $\GL_2$ \cite{serre-duke}.
This perspective of cohomology actually gives a natural explanation for the congruences for $\GL_2$ in \S \ref{sec:classical}. We explain how it naturally leads to considerations in modular representation theory.
From the short exact sequence $0 \ra W \overset{\cdot p}{\ra} W \ra W\otimes_{\Z_p} \F_p \ra 0$, we see that $H^*(Y(K_f),\cW)_{\fm}$ is nonzero if and only if $H^*(Y(K_f),\cW\otimes_{\Z_p}\F_p)_{\fm}$ is.
While $W(a,b) \otimes_{\Z_p} \Q_p$ is an irreducible $\GL_2(\Z_p)$-module, $W(a,b) \otimes_{\Z_p} \F_p$ is irreducible if and only if $a-b \leq p-1$, in which case, we let $F(a,b)\defeq W(a,b) \otimes_{\Z_p} \F_p$.
All (absolutely) irreducible $\GL_2(\Z_p)$-modules over $\F_p$ arise in this way, and $F(a,b) \cong F(c,d)$ if and only if $a-c = b-d \in (p-1)\Z$. 
Let $\cF(a,b)$ be the corresponding local system. 

Let us first revisit the congruences between Eisenstein series.
For any $a' > b$, the submodule of $W(a',b) \otimes_{\Z_p} \F_p$ generated by a (nonzero) highest weight vector is isomorphic to $F(a,b)$ where $a$ is the unique integer such that $0<a - b<p$ and $a \equiv a' \pmod{p-1}$.
This gives a map 
\[
H^*(Y(\GL_2(\widehat{\Z})),\cW(k-2,0) \otimes_{\Z_p} \F_p)_{\fm_{G_k}} \ra H^*(Y(\GL_2(\widehat{\Z})),\cW(k'-2,0)\otimes_{\Z_p} \F_p)_{\fm_{G_k}}
\]
where $2 < k < p+2$ and $k' > 2$ with $k' \equiv k \pmod{p-1}$.
It can be shown (for example by applying Hida's ordinary projector) that these maps are injective in each degree. 
This illustrates how modular representation theory can be used to produce infinite families of congruences between Hecke eigensystems.

The congruence \eqref{eqn:cuspcong} between cuspforms is simpler.
Recall that here $p=11$ and that $\bG = \GL_2$.
Then the fact that $\fm_\Delta$ is \emph{non-Eisenstein} implies that for all local systems $\cW$, $H^*(Y(K_f),\cW)_{\fm_\Delta}$ is zero unless $* = 1$.
First, Shapiro's lemma now implies that $H^*(Y(K_f),\cW)_{\fm_\Delta} \cong H^*(Y(\GL_2(\widehat{\Z})),\cW')_{\fm_\Delta}$, where $\cW'$ corresponds to the principal series representation $\Ind_{B(\F_p)}^{\GL_2(\F_p)} W$.
Second, the functor $W \mapsto H^1(Y(\GL_2(\widehat{\Z})),\cW)_{\fm_\Delta}$ is an exact functor from the category of finite $\Z_{11}[\![\GL_2(\Z_{11})]\!]$-modules to the category of finite $\Z_{11}$-modules.
Now $\Ind_{B(\F_{11})}^{\GL_2(\F_{11})} \bf{1}$ is naturally identified with the space of $\F_{11}$-valued functions on $\bP^1(\F_{11})$ and decomposes as $\Sym^{10} \F_{11}^2 \oplus \bf{1}$.
Then the injection
\[
H^*(Y(\GL_2(\widehat{\Z})),\cF(10,0))_{\fm_\Delta} \into H^*(Y(K_f),\F_p)_{\fm_\Delta}
\]
provides the desired congruence.
This example illustrates the important phenomenon of how the weight and level of modular forms can interact mod $p$.

With these representation theoretic arguments in mind, Ash, Stevens, and others have suggested that one should narrow the focus of Question \ref{ques:Serre1} to when $K_p$ is a maximal compact open subgroup and $\cW$ is an irreducible $\F$-local system.

\begin{ques}[The weight part of Serre's conjecture] \label{ques:Serre2}
Suppose that $K_p$ is a maximal compact open subgroup. Given a mod $p$ Hecke eigensystem $\fm$, for which irreducible $\F$-local systems $\cW$ on $Y(K_f)$ is $H^*(Y(K_f),\cW)_{\fm}$ nonzero?
\end{ques}

\noindent
While this is a substantial reduction since there are only finitely many such $\F$-local systems up to isomorphism, little is expected to be lost as we now explain.
The following proposition is immediate.

\begin{prop}\label{prop:JHnonzero}
If $H^*(Y(K_f),\cW)_{\fm}$ is nonzero, then $H^*(Y(K_f),\cF)_{\fm}$ is nonzero for some irreducible subquotient $\cF$ of $\cW \otimes_{\cO} \F$.
\end{prop}

\noindent The converse to Proposition \ref{prop:JHnonzero} holds in non-Eisenstein cases for $\bG$ a Weil restriction of $\GL_n$ if expected vanishing conjectures hold (see \S \ref{sec:vanish}).
These vanishing conjectures generalize the vanishing outside of degree $1$ for $\bG = \GL_2$.

Question \ref{ques:Serre2} turns out to be quite subtle.
Nonisomorphic irreducible $\F$-local systems may contain the same Hecke eigensystems, i.e.~not all congruences arise from modular representation theory.
For example, with $p=23$, both $H^*(Y(\GL_2(\widehat{\Z})),\cF(10,0))_{\fm_\Delta}$ and $H^*(Y(\GL_2(\widehat{\Z})),\cF(21,11))_{\fm_\Delta}$ are nonzero. 
If we write 
\[
\Delta(z) = \sum_{n=1}^\infty \tau(n) q^n,
\]
then $T_{\ell}$ acts on $H^*(Y(\GL_2(\widehat{\Z})),\cF(21,11))_{\fm_\Delta}$ by $\ell^{11}\tau(\ell)$ for all primes $\ell \neq 23$ by \eqref{eqn:matsushima} (see \S \ref{sec:classical}). This implies that
\begin{equation} \label{eqn:taucong}
\tau(n) \equiv \big(\tfrac{n}{23}\big)\tau(n) \pmod{23}
\end{equation}
for all $n$ coprime to $23$.
In other words, $23\mid \tau(n)$ if $\big(\tfrac{n}{23}\big) = -1$.

\section{An interlude on representation theory}\label{sec:reptheory}

\subsection{Serre weights} \label{sec:serrewt}

In order to explore Question \ref{ques:Serre2}, the natural first step is to ask for a classification of simple $\F[\![K_p]\!]$-modules.
If $K_p(1) \subset K_p$ is a normal finite index pro-$p$ subgroup and $W$ is a finite $\F[\![K_p]\!]$-module, then $W^{K_p(1)}$ is an $\F[\![K_p]\!]$-submodule of $W$ which is \emph{nonzero} (since the action of a $p$-group on a finite-dimensional $\F$-vector space must have a nonzero fixed vector).
Then the action of $K_p$ on a simple $\F[\![K_p]\!]$-module factors through the finite quotient $K_p/K_p(1)$, which can often be arranged to be a finite group of Lie type.
In this section, we discuss the (modular) representation theory of these groups.

Let $G$ be a connected reductive group over $\F_p$ which splits over $\F$. 
(We will eventually take $G$ to be the mod $p$ reduction of an integral model of $\bG$.)
An isomorphism class of a simple $\F[G(\F_p)]$-module is known as a \emph{Serre weight for} $G(\F_p)$.
Our goal now is to describe the (finite) set of Serre weights for $G(\F_p)$.

Let $B$ be an $\Fp$-rational Borel subgroup in $G$ with Levi subgroup $T$.
We denote by $W$ the Weyl group $N(T)/T$, which has a Bruhat partial order with a unique longest element $w_0$.
We write $X(T)$ for the character group of $T$, which has an action of $W$ and an induced action from the relative Frobenius $F$ acting on $T$.
This group has a subset 
\[
X_1(T) \defeq \{ \lambda \in X(T) : 0 \leq \langle \lambda,\alpha^\vee \rangle \leq p-1,\, \textrm{for all simple }\alpha \}
\]
which plays an important role in the modular representation theory of $G$.
We let 
\[
X^0(T) \defeq \{ \lambda \in X(T) : \langle \lambda,\alpha^\vee \rangle =0\, \textrm{ for all roots }\alpha \}.
\]

For any character $\lambda \in X(T)$, we can consider the algebraic induction $W(\lambda) \defeq \Ind_B^G w_0\lambda$ (also known as the dual Weyl module), which is nonzero if and only if $\lambda$ is dominant with respect to $B$.
We let $L(\lambda)$ denote the socle of $W(\lambda)$, which is the simple submodule generated by a nonzero highest weight vector.
Then we have the following result about Serre weights for $G(\F_p)$.

\begin{thm}
The map
\begin{align}
\frac{X_1(T)}{(F-1)X^0(T)} &\ra \{\textrm{Serre weights for } G(\F_p)\}\\
\lambda &\mapsto L(\lambda)(\F)|_{G(\F_p)}
\end{align}
is a bijection.
\end{thm}

\noindent We denote $L(\lambda)(\F)|_{G(\F_p)}$ by $F(\lambda)$.
(To avoid conflicts with the $F$-action on $X(T)$, we will write this action without parentheses.)

\subsection{Deligne--Lusztig representations}

While Question \ref{ques:Serre2} only involves $\F$-local systems, we will see that it is inextricably linked to $\cO$-torsion free local systems.
It is then natural to ask for a classification of irreducible $G(\F_p)$-representations in characteristic $0$.
We now recall such a classification, provided by work of Deligne and Lusztig \cite{DeligneLusztig}.

For an element $w \in W$, there exists $g_w \in G(\ovl{\F}_p)$ such that $g_w^{-1} F(g_w) \in N(T)(\ovl{\F}_p)$ represents $w$.
Then we let $T_w$ be the $F$-stable torus $g_w T g_w^{-1}$.
Let $\tld{W}$ denote the extended affine Weyl group which is the semidirect product $X(T) \rtimes W$.
For an element $\tld{w} = (\mu, w) \in \tld{W}$, we define $\theta_{\tld{w}}: T_w(\F_p) \ra E^\times$ to be the restriction of the character
\begin{align*}
T_w(\ovl{\F}_p) &\ra \ovl{\Q}_p^\times\\
g &\mapsto [\mu](g_w^{-1}g g_w)
\end{align*}
to $T_w(\F_p)$ (here, $[\mu]$ denotes the Teichm\"uller lift of $\mu$).

To a character $\theta_{\tld{w}}$ of a maximal rational torus $T_w(\F_p)$ of $G(\F_p)$, Deligne and Lusztig associate a virtual (Deligne--Lusztig) representation over $E$ which they denote $\epsilon_G \epsilon_{T_w} R_{T_w}^{\theta_{\tld{w}}}$.
We will instead denote this virtual representation by $R(\tld{w})$ and say that $\tld{w}$ is a presentation for $R(\tld{w})$.
The map $\tld{w} \mapsto R(\tld{w})$ is not injective---two elements map to the same virtual representation if and only if they lie in the same orbit of the action of $\tld{W}$ on itself given by
\[
(\nu,s)\cdot(\mu,w) = (s\mu+F\nu-swF(s)^{-1}(\nu),swF(s)^{-1}).
\]
The simplest case of the above construction occurs when $w$ is the identity. 
Then $T_w = T$, $\theta_{\tld{w}}$ is a character of $T(\F_p)$ and by inflation a character of $B(\F_p)$, and $R(\tld{w})$ is the principal series representation $\Ind_{B(\F_p)}^{G(\F_p)} \theta_{\tld{w}}$. 
Nonuniqueness of presentations can be seen from the existence of intertwiners between principal series representations.

The group $\tld{W}$ acts on $X(T)$ in the usual way---$W$ acts on $X(T)$ by group automorphisms and $X(T)$ acts on itself by translation. Let $m$ be a nonnegative integer and let $0 \in X(T)$ denote the trivial character.
We say that $\tld{w} \in \tld{W}$ is (lowest alcove) $m$-generic if $\langle \tld{w}(0),\alpha^\vee \rangle > m$ for all simple roots $\alpha$ and $\langle \tld{w}(0),\alpha^\vee \rangle < p-m$ for all roots $\alpha^\vee$.
We say that a Deligne--Lusztig representation $R$ is $m$-generic if $R = R(\tld{w})$ for some $m$-generic $\tld{w}$. 
An $m$-generic $\tld{w}$ or $R$ exists only if $m h < p$, where $h$ denotes the Coxeter number of $G$.
\cite[Proposition 10.10]{DeligneLusztig} implies that if $R$ is $0$-generic, then $R$ is in fact a genuine representation.

Let $\ovl{R}$ denote the semisimplification of the reduction of any $G(\F_p)$-stable $\cO$-lattice in a genuine $G(\F_p)$-representation $R$ over $E$ ($\ovl{R}$ does not depend on the choice of lattice).
In relation to Question \ref{ques:Serre2}, it is important to have an understanding of $\ovl{R}$ for Deligne--Lusztig representations $R$.
This is provided by Jantzen's formula for the reductions of Deligne--Lusztig representations in terms of virtual linear combinations of dual Weyl modules \cite{jantzenDL}.
If $R$ is sufficiently generic, then the Jordan--H\"older factors of $\ovl{R}$ admit the following description in terms of alcove geometry, which is in a sense independent of $p$.
For convenience, we assume that $G$ admits a \emph{twisting element} $\eta \in X(T)$, defined up to $X^0(T)$, which by definition has the property that $\langle \eta,\alpha^\vee\rangle = 1$ for all simple roots $\alpha$.
The existence of an $\eta$ can always be arranged by passing to a central extension of $G$ by $\bG_m$ (see \cite[Proposition 5.3.1(a)]{BG}).
We write $\cdot$ for the $p$-dot action so that $(\nu,w)\cdot \lambda = w(\lambda+\eta) - \eta + p\nu$. 
See \cite[\S 2]{MLM} for any unexplained notation below.

\begin{prop} \label{prop:JH} \cite[Proposition 2.3.6]{MLM}
Let $h$ be the Coxeter number of $G$.
If $\tld{w} \in \tld{W}$ is $2h$-generic, then the Jordan--H\"older factors of $\ovl{R}(\tld{w})$ are precisely the Serre weights of the form
\[
F(\pi^{-1}(\tld{w}_1)\cdot(\tld{w}\tld{w}_2^{-1}(0)-\eta))
\]
with $\tld{w}_1\in \tld{W}$ restricted and dominant, $\tld{w}_2 \in \tld{W}$ dominant, and $\tld{w}_1 \uparrow (\eta,w_0) \tld{w}_2$.
\end{prop}

\begin{rmk}
Of course, the description in Proposition \ref{prop:JH} does not depend on the choice of twisting element $\eta$ and could in fact be rephrased without any reference to $\eta$. 
\end{rmk}

\section{Relations to Galois representations}\label{sec:galrep}

In order to address Question \ref{ques:Serre2} (and to explain the congruence \eqref{eqn:taucong}), we introduce some conjectures and results concerning the relationship between cohomological automorphic forms and Galois representations.
We follow the approach in \cite{gross}, which seems to be more standard when $\bG$ is a general linear group or a unitary group.
For a more canonical approach to conjectures concerning Galois representations attached to cohomological automorphic forms, see \cite{BG}.

\subsection{Twisting element}

For a field $F$, let $G_F$ denote the absolute Galois group $\Gal(F^{\mathrm{sep}}/F)$ where we fix some separable closure $F^{\mathrm{sep}}$. 
Fix a maximal torus $T$ and Borel subgroup $B$ in $\bG_{/\ovl{\Q}}$.
We assume now that $\bG$ has a \emph{twisting element} $\eta$ which is by definition an element of $X(T)^{G_{\Q}}$ such that $\langle \eta,\alpha^\vee\rangle=1$ for all simple roots $\alpha$. 
If $\bG$ is $\GL_n$, we can take $\eta$ to be $(n-1,n-2,\ldots,1,0)$.
As before, a twisting element always exists if we replace $\bG$ by a central extension of $\bG$ by $\bG_m$.
The effect of this on the constructions and questions in \S \ref{sec:coh} is minimal. 
A twisting element is only unique up to $X^0(T)^{G_{\Q}}$.

\subsection{Satake parameters}\label{sec:satake}

Fix a prime $\ell$ and suppose that there is a reductive model $G$ over $\Z_\ell$ for $\bG$ such that $G(\Z_\ell) = K_\ell$.
Let $T\subset B \subset G$ be a maximal torus and Borel subgroup, respectively.
We have the Satake isomorphism (see e.g.~\cite[\S 4.2]{cartier})
\[
\cS: \Z[1/\ell][K_\ell \backslash \bG(\Q_\ell) / K_\ell] \risom \Z[1/\ell][T(\Z_\ell) \backslash T(\Q_\ell) / T(\Z_\ell)]^{W_s}
\]
normalized using the choice of twisting element $\eta$ as in \cite[Proposition 3.6]{gross-satake}, where $W_s$ is the Weyl group of the maximal split torus in $T$.

If $T$ is split, then $T(\Z_\ell) \backslash T(\Q_\ell) / T(\Z_\ell) \cong Y(T)$ (here $Y(T)$ denotes the cocharacter group of $T$) and $E$-valued characters of 
\[
\Z[1/\ell][K_\ell \backslash \bG(\Q_\ell) / K_\ell]\cong \Z[1/\ell][Y(T)]^W \cong \cO(\widehat{T}/\!\!/W)
\]
are in bijection with semisimple conjugacy classes in the dual group $\widehat{G}(E)$ for any coefficient field $E$ of characteristic not equal to $\ell$.
In general, $E$-valued characters $\chi$ of $\Z[1/\ell][K_\ell \backslash \bG(\Q_\ell) / K_\ell]$ are in bijection with semisimple conjugacy classes $C_\chi$ of $^L G(E)$, where $^L G$ denotes the Langlands dual of $G$ (see \cite[\S 16]{gross}).

\subsection{Conjectures on Galois representations associated to cohomological automorphic forms}

We fix a prime $p$ and a sufficiently large subfield $E \subset \ovl{\Q}_p$ of finite degree over $\Q_p$.

\begin{conj}\label{conj:galrep}
Let $V(\lambda)_E$ denote the irreducible representation of $\bG_{/E}$ of highest weight $\lambda$.
Suppose that $\fp \subset T_{\Q}^\Sigma \otimes_{\Q} E$ is a maximal ideal such that the $\fp$-torsion $H^*(Y(K_f),\cV(\lambda)_E) [\fp]$ is nonzero.
Then there exists a continuous homomorphism
\[
\rho: G_{\Q} \ra\, ^L\bG(E)
\]
such that 
\begin{enumerate}
\item the composition of $\rho$ with the projection $^L\bG(E) \ra G_{\Q}$ is the identity on $G_{\Q}$;
\item for $\ell \notin \Sigma$ and $\ell \neq p$, $\rho$ is unramified at $\ell$ and $\rho(\Frob_\ell)$ is in $C_\chi$ where $\chi$ is the character 
\[
\Z[1/\ell][K_\ell \backslash \bG(\Q_\ell) / K_\ell] \subset T_{\Q}^\Sigma \otimes_{\Q} E \ra (T_{\Q}^\Sigma \otimes_{\Q} E)/\fp \cong E
\] and $\Frob_\ell$ is an arithmetic Frobenius element at $\ell$; 
\item $\rho|_{G_{\Q_p}}$ is de Rham with Hodge--Tate cocharacter $\lambda+\eta$ (see \cite[\S 2.4]{BG}; in our normalization the cyclotomic character corresponds to the cocharacter $\id_{\bG_m}$) and is moreover crystalline if $p \notin \Sigma$; and
\item $\rho$ is odd in the sense of \cite[Conjecture 17.2(a)]{gross}.
\end{enumerate}
\end{conj}

There is an analogous conjecture with torsion coefficients. As before, let $\F$ denote the residue field of $E$.

\begin{conj}\label{conj:modp_galrep}
If $W$ is a finite $\F[\![K_p]\!]$-module, let $\cW$ denote the $\F$-local system on $Y(K_f)$.
Suppose that $\fm \subset T_{\cO}^\Sigma$ is a maximal ideal such that $H^*(Y(K_f),\cW)_{\fm}$ is nonzero.
Then there exists a continuous homomorphism
\[
\rhobar: G_{\Q} \ra \, ^L\bG(\F)
\]
such that 
\begin{enumerate}
\item the composition of $\rhobar$ with the projection $^L\bG(\F) \ra G_{\Q}$ is the identity on $G_{\Q}$;
\item for $\ell \notin \Sigma$ and $\ell \neq p$, $\rhobar$ is unramified at $\ell$ and $\rhobar(\Frob_\ell)$ is in $C_\chi$ where $\chi$ is the character \[
\Z[1/\ell][K_\ell \backslash \bG(\Q_\ell) / K_\ell] \ra T_{\Q}^\Sigma \otimes_{\Q}\F  \ra T_{\Q}^\Sigma \otimes_{\Q} \F/\fm \cong \F
\] and $\Frob_\ell$ is an arithmetic Frobenius element at $\ell$; and
\item $\rhobar$ is odd in the sense of \cite[Conjecture 17.2(a)]{gross}.
\end{enumerate}
\end{conj}

\begin{rmk}\label{rmk:galrep}
\begin{enumerate}
\item One also expects that $\rho$ satisfies a compatibility with the conjectural local Langlands correspondence at places in $\Sigma$.
\item Comparing the two conjectures, observe that there is no property at $p$ for $\rhobar$.
Such a property would be closely related to Questions \ref{ques:Serre1} and \ref{ques:Serre2}.
\item \label{item:cheb} It is clear that $\rho$ and $\rhobar$ determine $\fp$ and $\mathfrak{m}$, respectively.
On the other hand, the properties of $\rho$ described in Conjectures \ref{conj:galrep} and \ref{conj:modp_galrep} do not characterize $\rho$ in general.
When $\bG$ is a Weil restriction of $\GL_n$, the first two properties characterize the isomorphism class of the semisimplification of $\rho$ by the Chebotarev density theorem and the Brauer--Nesbitt theorem.
But even for tori, these properties do not characterize $\rho$ (see \cite[Remark 3.2.4]{BG}). 
\item The existence of torsion cohomology classes means that Conjecture \ref{conj:galrep} does not immediately imply Conjecture \ref{conj:modp_galrep}.
A version of Conjecture \ref{conj:modp_galrep} for all torsion coefficients implies Conjecture \ref{conj:galrep} (except for the third property) by taking a limit.
\end{enumerate}
\end{rmk}

There are two cases, relevant to what follows, when both conjectures are known: 
\begin{enumerate}
\item the Weil restriction of a definite unitary group relative to a CM extension of a totally real field not equal to $\Q$ \cite{kottwitz,HT,labesse,shin,CH}, and \\
\item the Weil restriction of $\GL_n$ over a CM field \cite{scholze}.
\end{enumerate}

In these cases, the attached Galois representations are determined up to semisimplification by Remark \ref{rmk:galrep}\eqref{item:cheb}.
We say that $\rhobar$ (or $\fm$) is non-Eisenstein if $\rhobar$ does not factor through a proper parabolic subgroup after any finite extension of $\F$. 

\subsection{Modular Serre weights}

Fix $p$ and $E$ as before, and let $\F$ be the residue field of $E$. 
Suppose from now on that $\bG_{/\Q_p}$ has an integral model $G_{/\Z_p}$. 
Having classified irreducible representations of $G(\Z_p)$ and introduced Galois representations, we now revisit Question \ref{ques:Serre2} through that lens. 
We let $K_p$ be $G(\Z_p)$. 
An irreducible $G(\Z_p)$-representation over $\F$ factors through the reduction map $G(\Z_p) \surj G(\F_p)$ (whose kernel is pro-$p$), i.e.~is the inflation of a Serre weight for $G(\F_p)$. 
For a Serre weight $\sigma$, let $\cF_\sigma$ denote the corresponding $\F$-local system on $Y(K_f)$. 

Fix an $\F$-valued Hecke eigensystem $\fm$. 
To understand Question \ref{ques:Serre2} is to understand the following (finite) set. 
\begin{defn}
Let $W(\fm)$ be the set of isomorphism classes of Serre weights $\sigma$ for $G(\F_p)$ for which $H^*(Y(K_f),\cF_\sigma)_\fm$ is nonzero. 
If $\sigma \in W(\fm)$, then we say that $\sigma$ is a \emph{modular Serre weight} (for $\fm$). 
\end{defn}

\noindent If there is a Galois representation $\rhobar: G_{\Q} \ra \, ^L\bG(\F)$ satisfying the properties in Conjecture \ref{conj:modp_galrep} for $\fm$, then we also write $W(\rhobar)$ for $W(\fm)$ and say that $\sigma$ is a modular Serre weight for $\rhobar$ when $\sigma \in W(\rhobar)$. 

\subsection{The case of $\GL_2$} \label{sec:GL2}

Fix $p$ and $E$ as before.
The first result on Conjecture \ref{conj:galrep} for nonabelian $\bG$ was work of Deligne \cite{deligne} in the case of $\bG=\GL_2$ (building on work of Eichler and Shimura).
In this case there is no torsion in cohomology, and so Conjecture \ref{conj:galrep} implies Conjecture \ref{conj:modp_galrep}.
Deligne constructed $\rho$ satisfying all the properties of Conjecture \ref{conj:galrep} except the third, which follows from subsequent work of Faltings on $p$-adic comparison theorems.
(We set $\eta = (1,0)$ here.)

Let $K_p = \GL_2(\Z_p)$.
Fix a mod $p$ Hecke eigensystem $\fm$. 
Assume that $\fm$ is non-Eisenstein, i.e.~that the attached Galois representation $\rhobar: G_{\Q} \ra \GL_2(\F)$ is absolutely irreducible. 
(The data of $\fm$ is equivalent to that of the isomorphism class of $\rhobar$ by Remark \ref{rmk:galrep}\eqref{item:cheb}.)
Since cohomology groups outside degree $1$ do not admit non-Eisenstein mod $p$ Hecke eigensystems, the functor $W \ra H^1(Y(K_f),\cW)_{\fm}$ is exact, as explained in \S \ref{sec:congruences}.
In particular, Question \ref{ques:Serre1} reduces to Question \ref{ques:Serre2}, namely an investigation of $W(\rhobar)$. 


If $H^1(Y(K_f),\cF(a,b))_\fm$ is nonzero, then so is $H^1(Y(K_f),\cW(a,b))_\fm$.
A necessary condition for $F(a,b)$ to be in $W(\rhobar)$ is that $\rhobar$ is the reduction of a representation $\rho: G_{\Q} \ra \GL_2(E)$ which is unramified outside of $\Sigma$ and $p$ and is crystalline at $p$ of Hodge--Tate weights $a+1$ and $b$.
(Since $\rhobar$ is irreducible, a $G_{\Q}$-invariant $\cO$-lattice in $\rho$ is unique up to scaling.)
In particular, the restriction $\rhobar|_{G_{\Q_p}}$ is the reduction of (a lattice in) a crystalline representation $\rho_p: G_{\Q_p} \ra \GL_2(E)$ of Hodge--Tate weights $a+1$ and $b$.
The following result, known as the weight part of Serre's conjecture, is a local-global principle (i.e.~$W(\rhobar)$ only depends on $\rhobar|_{G_{\Q_p}}$) that asserts that the necessary condition is in fact sufficient.

\begin{thm}[\cite{gross,edixhoven,CV}]\label{thm:SWC_GL2}
Suppose that $\fm$ is non-Eisenstein.
Then $F(a,b)\in W(\rhobar)$ if and only if $\rhobar|_{G_{\Q_p}}$ is the reduction of a crystalline representation of Hodge--Tate weights $a+1$ and $b$. 
\end{thm}

\begin{rmk}\label{rmk:SWC_GL2}
\begin{enumerate}
\item The above formulation is slightly different, albeit equivalent, from Serre's original formulation in \cite{serre-duke}. In \emph{loc.cit.}, the recipe for $W(\rhobar)$ is completely explicit in terms of the ``inertial weights" of $\rhobar|_{G_{\Q_p}}$ when $\rhobar|_{G_{\Q_p}}$ is semisimple, with additional modfications in terms of ramification properties of an extension class in general.
In particular, in the semisimple case, there is a simple  combinatorial formula for $a$ and $b$ solely in terms of the inertial weights.
Early generalizations of Serre's conjecture beyond $\GL_{2/\Q}$, e.g.~\cite[Conjecture 3.1]{ADP}, involved similar formulas (see \cite[\S 7]{GHS} for more recent formulas for a larger list of Serre weights). 
On the other hand, while \cite{BDJ} also contains formulas \`a la Serre, it emphasizes the above ``crystalline lifts" perspective. 
\item Theorem \ref{thm:SWC_GL2} was generalized to (the Weil restriction of) the unit groups in quaternion algebras over totally real fields split at no more than one archimedean place and definite unitary groups over a totally real field when $p>2$ (and under mild additional hypotheses) in a series of works by Gee, Newton, Kisin, Liu, and Savitt \cite{newton13,gee-kisin,GLS}.
Some of these build on earlier work of Gee that introduced the Taylor--Wiles method to produce proofs rather different from the original proofs of Theorem \ref{thm:SWC_GL2}.
\end{enumerate}
\end{rmk}

We now revisit the mod $23$ congruence for the Ramanujan Delta function.
Let $\rhobar: G_{\Q} \ra \GL_2(\F_{23})$ be the associated Galois representation.
In this case, $\rhobar|_{I_{\Q_{23}}} \cong \omega^{11} \oplus \bf{1}$, where $\omega$ denotes the reduction of the $23$-adic cyclotomic character $\chi$ and $I_{\Q_{23}}\subset G_{\Q_{23}}$ denotes the inertial subgroup.
Then $\rhobar|_{G_{\Q_{23}}}$ is the reduction of both $\chi^{11} \oplus \bf{1}$ and $\chi^{11} \oplus \chi^{22}$ (up to unramified twists).
Theorem \ref{thm:SWC_GL2} implies that $\{F(10,0),F(21,11)\} \subset W(\rhobar)$.
In fact, this is an equality.

The behavior of the Ramanujan Delta function modulo $23$ illustrates a rare phenomenon.
For example, if we instead take $p = 19$ and $\rhobar: G_{\Q} \ra \GL_2(\F_{19})$ the corresponding representation, then $\rhobar|_{I_{\Q_{19}}}$ is a \emph{nontrivial} extension of $\bf{1}$ by $\omega^{11}$. 
Moreover, $\rhobar|_{G_{\Q_{19}}}$ is the reduction of a crystalline representation which after restriction to inertia is a nontrivial extension of $\bf{1}$ by $\chi^{11}$.
However, any nontrivial extension of $\chi^{18}$ by any unramified twist of $\chi^{11}$ is \emph{not} crystalline.
In fact, we have $\{F(10,0)\} = W(\rhobar)$ in this case.
One expects more generally that $W(\rhobar)$ is larger when $\rhobar|_{G_{\Q_p}}$ is semisimple.
Indeed, if $\rhobar|_{G_{\Q_p}}$ is the reduction of an $\cO$-lattice in $r: G_{\Q_p} \ra \GL_2(E)$, then there is another $\cO$-lattice (possibly after enlarging $E$) whose reduction is the semisimplification of $\rhobar|_{G_{\Q_p}}$.

One approach to Theorem \ref{thm:SWC_GL2} is as follows.
If $\rhobar|_{G_{\Q_p}}$ admits a local crystalline lift of Hodge--Tate weights $a+1$ and $b$, show that $\rhobar$ admits a global lift $\rho$ which is crystalline (at $p$) of Hodge--Tate weights $a+1$ and $b$.
Then show that $\rho$ comes from a modular form as predicted by the Fontaine--Mazur conjecture \cite{fontaine-mazur}.
These steps can be executed using a combination of tools in the Taylor--Wiles method independently discovered by Khare--Wintenberger and Gee \cite{KW_Serre_2,gee-annalen}.

\subsection{Vanishing conjectures for cohomology}\label{sec:vanish}

In the previous section, we were in the advantageous situation where $H^*(Y(K_f),\cW)_{\fm}$ vanished outside of degree one for all $\cO$-local systems $\cW$ when $\fm$ is non-Eisenstein. 
While this is not true in general, this property does nevertherless admit a conjectural generalization.
Let $d_Y$ be the dimension of $Y(K_f)$.
Define the integers $\ell_0\defeq \rk \bG(\R) - \rk A_\infty^\circ K_\infty^\circ$ and $q_0 = \frac{1}{2}(d_Y-\ell_0)$.
(The group $\bG$ admits discrete series if and only of $\ell_0 = 0$.)

\begin{conj}[\cite{CG}] \label{conj:vanish}
Suppose that $\fm$ is non-Eisenstein.
If $H^i(Y(K_f),\cW)_{\fm} \neq 0$, then $i \in [q_0,q_0+\ell_0]$.
\end{conj}

\noindent \cite{GN} shows that if Conjecture \ref{conj:vanish} holds, then so does the converse to Proposition \ref{prop:JHnonzero} when $\bG$ is a Weil restriction of $\GL_n$, so that, as in our discussion above, it suffices to analyze Question \ref{ques:Serre1} for irreducible mod $p$ local systems (note that this is far from obvious in general, as it is \emph{a priori} possible for the cohomology complex of irreducible local system to cancel each other out when they spread to several cohomological degrees). More seriously, Conjecture \ref{conj:vanish} plays a prominent role in the Taylor-Wiles patching method, which is the main tool to attack Question \ref{ques:Serre1} in general, cf.~ Remark \ref{rmk:patch} below.

Unfortunately, there are few cases where Conjecture \ref{conj:vanish} is known.
One trivial case is that of groups which are anisotropic modulo their center, e.g.~definite unitary groups, when $Y(K_f)$ is a finite set of points.
\cite{CS} have shown that for certain unitary groups ($\ell_0 = 0$) Conjecture \ref{conj:vanish} holds under some additional hypotheses.
The case of the Weil restriction of $\GL_n$ (for $n>1$) over a number field $F$ (even CM fields) is open beyond the case $n=2$ and $F$ either totally real or imaginary quadratic.

\section{Conjectures and results on the weight part of Serre's conjecture}\label{sec:conj}

\subsection{Taylor--Wiles patching}\label{sec:TW}

Suppose for the moment that we are in a context where $\ell_0= 0$ and Conjectures \ref{conj:galrep} and \ref{conj:vanish} hold (e.g., definite unitary groups). 
We fix $p$ and $E$ as before. 
We assume that a reductive integral model $G_{/\Z_p}$ of $\bG_{/\Q_p}$ exists and continue to let $K_p=G(\Z_p)$. 
If $F(\lambda) \in W(\rhobar)$, then $H^*(Y(K_f),\cW(\lambda))_\fm$ is nonzero (where $W(\lambda)$ is an $\cO$-lattice in the irreducible algebraic representation $V(\lambda)$).
In particular, $\rhobar$ (attached to $\fm$) is the reduction of a crystalline representation of Hodge--Tate cocharacter $\lambda+\eta$.
In light of Theorem \ref{thm:SWC_GL2}, it is tempting to guess that the converse holds.
However, counterexamples to this have been found for definite unitary groups in three variables \cite[Proposition 7.18]{LLLM}.
The reason is that in contrast to the case of $\GL_2$, $W(\lambda)$ may include many Jordan--H\"older factors other than $F(\lambda)$.
In fact, $F(\lambda)$ as a $\cO[\![K_p]\!]$-module often does not lift to an $\cO$-torsion-free module, which makes it more difficult to use the (expected) $p$-adic Hodge theoretic properties of Galois representations attached to automorphic forms.
However, $F(\lambda)$ does lift \emph{virtually}.

Suppose that $[F(\lambda)] = \sum_W c_W [W]$ in the Grothendieck group of $\F[\![K_p]\!]$-modules, where each $W$ in the sum lifts to characteristic $0$ (for example, one can take $W$ running over the reductions of various $\cO[G(\F_p)]$-modules using \cite[Theorem 33]{serre-book}).
Then exactness gives us 
\[
[H^*(Y(K_f),\cF(\lambda))_{\fm}] = \sum_W c_W [H^*(Y(K_f),\cW)_{\fm}].
\]
Since the $\F$-vector spaces on the right hand side lift, their dimensions can in principle be computed in characteristic zero.
However, they are of a global nature, and thus difficult to access.
In contrast, we still expect that $W(\rhobar)$ depends only on $\rhobar|_{G_{\Q_p}}$.
The Taylor--Wiles method ``patches" together cohomology functors (or rather, the total cohomology complex computing) $H^*(Y(K_f),\cF_-)_{\fm}$ (for varying $K_f$) to obtain a functor $M_\infty(-)$ that plausibly depends (roughly speaking) only on $\rhobar|_{G_{\Q_p}}$.
Moreover, a control theorem guarantees that $M_\infty(\sigma)$ is nonzero if and only if $H^*(Y(K_f),\cF_\sigma)_{\fm}$ is.

For a Galois representation $\rbar: G_{\Q_p}\ra \, ^L \bG(\F)$, let $R_{\rbar}$ denote the (framed) deformation ring parametrizing lifts $r: G_{\Q_p}\ra \, ^L \bG(R)$ of $\rbar$ for complete local Noetherian $\cO$-algebras $R$ with residue field $\F$.
Building on work of Kisin \cite{KisinPSS} in the case when $\bG_{/\Q_p}$ is a Weil restriction of $\GL_n$, Balaji \cite{balaji} in particular defined a family of (reduced) semistable deformation rings $R^{\lambda+\eta,\tau}_{\rbar}$ whose $\ovl{\Q}_p$-points correspond to potentially semistable Galois representations of Hodge--Tate cocharacter $\lambda+\eta$ and Galois type $\tau$.
In certain contexts where (enough of) an inertial local Langlands correspondence is known, one can define a finite dimensional locally algebraic $E[K_p]$-module $\sigma(\lambda,\tau) \defeq V(\lambda) \otimes_E \sigma(\tau)$.
For example, when $\tau$ is tame, $\sigma(\tau)$ can be taken to be a certain combinatorially defined Deligne--Lusztig representation.
For a ring $A$, let $A\textrm{-mod}^{\textrm{fg}}$ denote the full subcategory of $A$-mod of finitely generated $A$-modules.

\begin{ax}\label{ax:patch}
There is an exact functor $M_\infty(-): \cO[\![K_p]\!]\textrm{-mod}^{\textrm{fg}} \ra R_{\rhobar|_{G_{\Q_p}}}\textrm{-mod}^{\textrm{fg}}$ with the following properties.
\begin{enumerate}
\item \label{item:CM} For a Serre weight $\sigma$, $M_\infty(\sigma)$ is a maximal Cohen--Macaulay module on $R^{\lambda+\eta,\tau}_{\rbar} \otimes_{\cO} \F$ for some $\lambda$ and $\tau$.
\item \label{item:control} $M_\infty(\sigma)$ is nonzero if and only if $H^*(Y(K_f),\cF_\sigma)_{\fm}$ is nonzero.
\item \label{item:LGC} If $\sigma(\lambda,\tau)^\circ$ is an $\cO$-lattice in $\sigma(\lambda,\tau)$, then $M_\infty(\sigma(\lambda,\tau)^\circ)$ is a maximal Cohen--Macaulay $R^{\lambda+\eta,\tau}_{\rbar}$-module of generic rank at most $1$.
\item \label{item:fullsupport} If $\sigma(\lambda,\tau)^\circ$ is an $\cO$-lattice in $\sigma(\lambda,\tau)$, then $M_\infty(\sigma(\lambda,\tau)^\circ)[1/p]$ is a generically free $R^{\lambda+\eta,\tau}_{\rbar}[1/p]$-module of rank $1$.
\end{enumerate}
\end{ax}

\begin{rmk}\label{rmk:patch}
\begin{enumerate}
\item It may be impossible to arrange for the stringent rank conditions in Axiom \ref{ax:patch}, but the ranks can still be controlled and many arguments below can be successfully modified.
\item If Conjecture \ref{conj:vanish} holds, then it is often possible to use the Taylor--Wiles method to construct a functor $M_\infty(-)$ satisfying the first three items of Axiom \ref{ax:patch} with $\rbar \defeq \rhobar|_{G_{\Q_p}}$ \cite{CEGGPS,GN}, at least after adding formal variables to $R^{\lambda+\eta,\tau}_{\rbar}$ which we ignore. 
In general, the total (localized) cohomology complex may have non-vanishing cohomology groups in several degrees. The role of Conjecture \ref{conj:vanish} is to guarantee that after Taylor-Wiles patching, these complexes becomes concentrated into a single cohomological degree, i.e.~they turn into usual modules. This concentration effect relies on the ``numerical coincidence'' that powers the Taylor-Wiles method.
\item \label{item:domaincase} It seems to be difficult to guarantee the last item (for all choices of $(\lambda,\tau)$) except when $\bG = \GL_2$ \cite{kisin-fontaine-mazur}. 
Indeed, it is essentially equivalent to a modularity lifting result with very general $p$-adic Hodge theoretic hypotheses. 
However, there are some instances of specific $(\lambda,\tau)$ where the final item follows from the third: when $R^{\lambda+\eta,\tau}_{\rbar}$ is zero and when $R^{\lambda+\eta,\tau}_{\rbar}$ is a domain and $M_\infty(\sigma(\lambda,\tau)^\circ)$ is nonzero for some $\cO$-lattice $\sigma(\lambda,\tau)^\circ\subset \sigma(\lambda,\tau)$.
\end{enumerate}
\end{rmk}

Admitting Axiom \ref{ax:patch} for the moment, there is the following strategy for determining $W(\rhobar)$.
Let $d = \dim_{E} G_{/E}+\dim_{E} G_{/E}/B_{/E}$ where $B_{/E} \subset G_{/E}$ is a Borel subgroup.
By a result of Kisin (see Theorem \ref{thm:local_properties}), $d$ is the dimension of $R^{\lambda+\eta,\tau}_{\rbar}$ over $\cO$ for any $\lambda$ and $\tau$. 
For each Serre weight $\sigma$, we write a presentation
\begin{equation}\label{eqn:present}
[\sigma] = \sum_{(\lambda,\tau)} c_{\lambda,\tau} [{\ovl{\sigma(\lambda,\tau)}}]
\end{equation}
in the Grothendieck group of $\F[\![K_p]\!]$-modules. 
Then Axiom \ref{ax:patch} guarantees that 
\begin{equation}\label{eqn:BM}
Z(M_\infty(\sigma)) = \sum_{(\lambda,\tau)} c_{\lambda,\tau} Z(R^{\lambda+\eta,\tau}_{\rbar} \otimes_{\cO} \F),
\end{equation}
where $Z(-)$ corresponds to taking the $d$-dimensional support cycle of the $R_{\rbar}$-module (note that all modules involved have support of dimension $\leq d$).
Since the right hand side of \eqref{eqn:BM} only depends on $\rbar \defeq \rhobar|_{G_{\Q_p}}$, so does the left hand side.
In particular, this implies $W(\rhobar)$, being the set of $\sigma$ with $Z(M_\infty(\sigma))\neq 0$ by the Cohen--Macaulay property, depends only on $\rbar$ as expected.

Another feature of the situation is that there are many possible choices of expressions (\ref{eqn:present}), even if we restrict to the case when $c_{\lambda,\tau} =0$ for all $\lambda \neq 0$. 
Since the left hand side of (\ref{eqn:BM}) involves only $\sigma$, we get the surprising conclusion that the right hand side must be independent of the choice of expressions  (\ref{eqn:present}). In other words, there are many non-trivial relations between cycles of special fibers of potentially crystalline deformation rings.
We summarize the above arguments as the following conjecture.

\begin{conj}\label{conj:GHS}
\begin{enumerate}
\item \label{item:BM} \cite{BM,EG} The left hand side of \eqref{eqn:BM} is independent of the presentation in \eqref{eqn:present}.
\item \label{item:BMSW} \cite{GHS} For any presentation as in \eqref{eqn:present}, $\sigma \in W(\rhobar)$ if and only if the right hand side of \eqref{eqn:BM} is nonzero.
\end{enumerate}
\end{conj}

\begin{rmk}
Conjecture \ref{conj:GHS}\eqref{item:BM} is purely local in the sense that it only involves $\bG_{/\Q_p}$ and $\rbar$ while \eqref{item:BMSW} is global since it involves $\rhobar$.
However, both follow from Axiom \ref{ax:patch}. 
\end{rmk}

\subsection{Herzig's recipe}
The complexity of Conjecture \ref{conj:GHS} suggests that Question \ref{ques:Serre2} does not admit a simple answer.
Indeed, the case of $\GL_2$ is perhaps misleading because of the simplicity of the geometry and representation theory involved.
However, the class of semisimple representations $G_\Q \ra \, ^L \bG(\F)$ admits an essentially combinatorial classification, and so one could ask for a combinatorial description of $W(\rhobar)$ when $\rhobar|_{G_{\Q_p}}$ is semisimple which generalizes Serre's recipe (see Remark \ref{rmk:SWC_GL2}).
Herzig's recipe gives such a (conjectural) description.

We assume in this section that $\bG_{/\Q_p}$ is unramified, with reductive integral model $G_{/\Z_p}$.
We let $K_p$ be $G(\Z_p)$.
As before, we fix a twisting element $\eta$ of $\bG$.
A \emph{regular} Serre weight is a Serre weight $F(\lambda)$ with $0 \leq \langle \lambda,\alpha^\vee \rangle< p-1$ for all simple roots $\alpha$.
We define an involution $\cR$ on the set of regular Serre weights $F(\lambda)$ by the formula
\[
\cR(F(\lambda)) = F((-w_0(\eta),w_0)\cdot \lambda).
\]

If $\rbar: G_{\Q_p} \ra \, ^L\bG(\F)$ is semisimple, then the restriction $\rbar|_{I_{\Q_p}}$ to the inertial subgroup factors through a torus. 
Its Teichm\"uller lift $[\rbar|_{I_{\Q_p}}]: G_{\Q_p} \ra \, ^L\bG(E)$ is then a tame inertial type.
Recall that to a tame inertial type $\tau$ (defined over $E$), one can attach through a tame inertial local Langlands a Deligne--Lusztig representation $\sigma(\tau)$ of $G(\F_p)$ (defined over $E$).
Then $K_p$ acts on $\sigma(\tau)$ by inflation, and we let $\ovl{\sigma}(\tau)$ denote the semisimplification of the reduction of any $K_p$-stable $\cO$-lattice in $\sigma(\tau)$.
We have the following conjecture.

\begin{conj}[\cite{herzig-duke,GHS}]\label{conj:herzig}
The subset of regular Serre weights in $W(\rhobar)$ is $\cR(\JH(\ovl{\sigma}([\rhobar|_{I_{\Q_p}}])))$.
\end{conj}

\noindent Conjectures \ref{conj:GHS} and \ref{conj:herzig} are of a rather different nature.
For one thing, Conjecture \ref{conj:herzig} only applies to $\rhobar$ which are semisimple locally at $p$, which is when one expects $W(\rhobar)$ is largest. 
However, Conjecture \ref{conj:herzig} is rather more explicit than Conjecture \ref{conj:GHS} when combined with Proposition \ref{prop:JH}.

\subsection{Results on the weight part of Serre's conjecture} \label{sec:result}

Conjecture \ref{conj:herzig} and a weakened version of Conjecture \ref{conj:GHS} is known when $\bG$ is $\GL_2$, the Weil restriction of the unit group in a quaternion algebra which is indefinite at no more than one archimedean place, or the Weil restriction of a definite unitary group in two variables under mild hypotheses (see Remark \ref{rmk:SWC_GL2}).
Similar results \cite{LLLM,LLLM2,wild} are known for the Weil restrictions of definite unitary groups in three variables that are unramified at $p$ under an additional \emph{genericity} hypothesis.

Suppose now that $\bG_{/E}$ is a product of $\GL_n$ over a finite set $\cJ$.
For $\tld{w} \in \tld{W}$, $\tld{w}(0)$ is a tuple of elements in $\Z^n$ indexed by $\cJ$.
We write $\tld{w}(0)_j \in \Z^n$ for the element corresponding to $j\in \cJ$.
We say that a semisimple $\rbar: G_{\Q_p} \ra \, ^L\bG(\F)$ is sufficiently generic at $p$ if ${\sigma}([\rbar|_{I_{\Q_p}}]) = R(\tld{w})$ where $0 \leq \langle \tld{w}(0),\alpha^\vee\rangle \leq p$ for all simple roots $\alpha$ and $p \nmid P(\tld{w}(0)_j)$ for an implicit polynomial $P \in \Z[X_1,\ldots,X_n]$ which depends only on $n$ (and not on $p$ or $j$).
If $\rbar$ is not semisimple, we say that it is sufficiently generic if its semisimplification is.
In a precise sense, most $\rbar$'s are sufficiently generic as $p \ra \infty$.

\begin{thm}\label{thm:main}
\begin{enumerate} 
\item \label{item:genBM} \cite[Corollary 8.5.2]{MLM} Suppose that $\bG_{/\Q_p}$ is an unramified Weil restriction of $\GL_n$.
Then Conjecture \ref{conj:GHS}\eqref{item:BM}, restricted to presentations coming from Deligne--Lusztig representations, holds for sufficiently generic $\rbar$.
\item \label{item:unitSW} \cite[Theorem 9.1.6]{MLM} Suppose that $\bG$ is the Weil restriction of a definite unitary (over a nontrivial totally real extension of $\Q$) and that $\rhobar|_{G_{\Q_p}}$ is sufficiently generic and semisimple.
Then under mild Taylor--Wiles hypotheses, Conjectures \ref{conj:GHS}\eqref{item:BMSW} (for the restricted presentations in the previous item) and \ref{conj:herzig} hold.
\item \label{item:GLnSW} \cite{LL} Suppose that $\bG$ is the Weil restriction of $\GL_n$ over a CM field and that $\rhobar|_{G_{\Q_p}}$ is sufficiently generic and semisimple.
Suppose moreover that $H^*(Y(K_f),\cW)_\fm \otimes_{\cO} E$ is nonzero for some $\cW$ corresponding to a Deligne--Lusztig representation.
Then under mild Taylor--Wiles hypotheses, Conjectures \ref{conj:GHS}\eqref{item:BMSW} (for the restricted presentations in the previous item) and \ref{conj:herzig} hold.
\end{enumerate}
\end{thm}

A critical tool to prove Theorem \ref{thm:main} is the following. 

\begin{thm}\label{thm:domain}
Suppose that $\rbar: G_{\Qp} \ra \, ^L \bG(\F)$ is semisimple. 
\begin{enumerate} 
\item \cite[Theorem 7.3.2(2)]{MLM} If $\tau$ is a sufficiently generic tame inertial type, then $R_{\rbar}^{\eta,\tau}$ is an integral domain (if it is nonzero). 
\item \label{item:pdpatch} \cite[Proposition 6.2.7]{MLM} There exists a functor $M_\infty$ (up to adding formal variables) satisfying Axiom \ref{ax:patch}\eqref{item:CM}-\eqref{item:LGC}, with $\bG$ a Weil restriction of a definite unitary group, such that $M_\infty(\sigma(\tau)^\circ)$ is nonzero if $R_{\rbar}^{\eta,\tau}$ is nonzero for a sufficiently generic tame inertial type $\tau$.
\end{enumerate}
\end{thm}

\begin{rmk}
\begin{enumerate}
\item Theorem \ref{thm:domain}\eqref{item:pdpatch} combines the modularity of obvious weights \cite{LLL} and the coherence conjecture for local models of Shimura varieties \cite{Zhu_coherence}. 
\item In fact, Theorem \ref{thm:domain} also holds for any $\lambda+\eta$ with $\lambda$ dominant, though the implicit polynomial defining genericity depends on $\lambda$. See Theorem \ref{thm:main_def_ring}.
\end{enumerate}
\end{rmk}

As alluded to in Remark \ref{rmk:patch}\eqref{item:domaincase}, Theorem \ref{thm:domain} implies the existence of a functor $M_\infty$ satisfying Axiom \ref{ax:patch}\eqref{item:CM}-\eqref{item:fullsupport} restricted to cases when $\lambda = 0$ and $\tau$ is a generic tame inertial type.
The argument from \S \ref{sec:TW} shows that Theorem \ref{thm:domain} implies Theorem \ref{thm:main}\eqref{item:genBM} for sufficiently generic \emph{semisimple} $\rbar$. 
The nonsemisimple case follows from a simple argument using the global geometry of the Emerton--Gee stack relying on the fact that there is a semisimple $\rbar$ on every component (see Remark \ref{rmk:BM}(\ref{rmk:suffmany})). 
Moreover, Theorem \ref{thm:main}\eqref{item:unitSW} follows from Axiom \ref{ax:patch}\eqref{item:control}. 

\begin{rmk}
The final part of Theorem \ref{thm:main} is a bit more subtle.
When $\ell_0=0$, one expects $H^*(Y(K_f),\cW)_\fm$ to have little torsion. In contrast when $\ell_0>0$, one expects cohomology to be dominated by torsion and characteristic $0$ classes to be rare. 
This means that the lifting hypothesis, i.e.~that $H^*(Y(K_f),\cW)_\fm \otimes_{\cO} E$ is nonzero in Theorem \ref{thm:main}\eqref{item:GLnSW}, is quite restrictive. 
This condition could be removed if one knew Conjecture \ref{conj:vanish} (see Remark \ref{rmk:patch}).
In lieu of this, the lifting hypothesis can be used to make an argument with Euler characteristics of the functor $M_\infty$, whose image is a priori an object in $D^b(R_{\rhobar|_{G_{\Q_p}}}\textrm{-mod}^{\textrm{fg}})$, adopting Taylor's Ihara avoidance trick to this setting \cite{10author}.
\end{rmk}

\section{Local models for potentially crystalline deformation rings}\label{sec:LM}

The heart of the proof of Theorem \ref{thm:main} reduces, via the Taylor-Wiles method, to understanding the support of the patched modules $M_\infty(\sigma)$ in Axiom \ref{ax:patch}, and ultimately to geometric properties of the potentially crystalline deformation rings $R^{\lambda,\tau}_{\rbar}$. We achieve this by introducing and analyzing certain (finite type) group-theoretic moduli spaces which algebraize these deformation rings.

Fix a prime $p$. 
In this section, we will restrict to the case $\bG_{/\Qp}=\mathrm{Res}_{K/\Qp}\GL_n$ for an unramified extension $K=\bQ_{p^f}$ of $\bQ_p$. In particular, we have a reductive integral model $G_{/\Zp}=\mathrm{Res}_{\cO_K/\Zp}\GL_n$ of $\bG$. 
Recall from \S \ref{sec:integral} that $E$ is a sufficiently large finite extension of $\Q_p$ with ring of integers $\cO$, uniformizer $\varpi$, and residue field $\F$. 
 
\subsection{Potentially crystalline deformation rings}
Let $\rbar: G_K\to \GL_n(\F)$ be a mod $p$ local Galois representation. Recall that we have $R_{\rbar}$ the (framed) deformation ring that classifies lifts of $\rbar$, and the $\Qpbar$-points of $R_{\rbar}$ correspond to $p$-adic Galois representations of $G_K$ (lifting $\rbar$). Given a Hodge--Tate cocharacter $\lambda$ and inertial type $\tau$, 
one has the potentially crystalline deformation ring $R_{\rbar}^{\lambda,\tau}$ which is characterized as the unique reduced $p$-flat quotient of $R_{\rbar}$
whose $\ovl{\Q}_p$-points correspond to lifts $r: G_K \ra \GL_n(\ovl{\Q}_p)$ which are potentially crystalline of type $(\lambda,\tau)$ (i.e.~the Hodge--Tate weights of $r$ are given by $\lambda$ and $\mathrm{WD}(r)$ induces the inertial type $\tau$).
In the setting of $\GL_n$, these rings were first constructed by Kisin, who also established their basic properties \cite{KisinPSS}. 
\begin{thm}[Kisin]\label{thm:local_properties}
\begin{enumerate}
\item $R_{\rbar}^{\lambda,\tau}[\frac{1}{p}]$ is regular.
\item $\dim_{\cO} R_{\rbar}^{\lambda,\tau}=\dim_{E} G_{/E}+\dim_{E} G_{/E}/P_{\lambda/E}$, 
where $P_{\lambda}$ is the parabolic subgroup corresponding to $\lambda$. In particular $\dim_{\cO}R^{\lambda,\tau}$ is a constant $d$ as $\lambda$ varies over regular dominant cocharacters.
\end{enumerate}
\end{thm}

When $\lambda$ is \emph{regular dominant}, the rings $R_{\rbar}^{\lambda,\tau}$ play a pivotal role in the Taylor-Wiles method: they act on patched spaces of automorphic forms $M_{\infty}(\sigma(\lambda-\eta,\tau))$, which govern questions about modularity and congruences (cf.~Axiom \ref{ax:patch}\eqref{item:LGC}). Even better, they are maximal Cohen--Macaulay modules, and hence must be supported on a union of irreducible components.

For global applications, it is essential to understand global properties of $R_{\rbar}^{\lambda,\tau}$ such as irreducibility. This turns out to be a notoriously difficult problem. There are roughly two reasons for this. 
\begin{itemize}
\item Outside some special cases, $R_{\rbar}^{\lambda,\tau}$, being characterized by its $\Qpbar$-points, has no known moduli interpretation. This is related to the fact that integral $p$-adic Hodge theory is much less well understood than rational $p$-adic Hodge theory.
\item The internal structure of $R_{\rbar}^{\lambda,\tau}$ is intrinsically complicated in general. Thus, one can not expect to have completely explicit descriptions for all $\lambda$ and $\tau$.
\end{itemize}
The second point is best illustrated by the Breuil--M\'ezard conjecture, which quantifies the complexity of the special fibers of $R_{\rbar}^{\lambda+\eta,\tau}$ as $\lambda$ and $\tau$ vary in terms of the mod $p$ representation theory of $\GL_n(\cO_K)$ (we shift from $\lambda$ to $\lambda+\eta$ for the rest of this subsection to be consistent with \S \ref{sec:TW}).
We let $Z(R_{\rbar}^{\lambda+\eta, \tau}/\varpi)$ denote the $d$-dimensional cycle of $R_{\rbar}^{\lambda+\eta, \tau}/\varpi$, which counts the irreducible components with appropriate multiplicities. 
For a $\GL_n(\cO_K)$-representation $V$ over $E$, recall that $\ovl{V}$ denotes the $\GL_n(\cO_K)$-representation over $\F$ which is the semisimplification of the reduction modulo $\varpi$ of any $\GL_n(\cO_K)$-stable $\cO$-lattice in $V$. 
The following is a reformulation of Conjecture \ref{conj:GHS}\eqref{item:BM}.

\begin{conj}[Breuil--M\'ezard, Emerton-Gee]  \label{conj:introBM} There exist $d$-dimensional cycles $\cZ_{\sigma}(\rbar)$ in $\Spec R_{\rbar}/\varpi$ for each irreducible $\GL_n(\cO_K)$-representation $\sigma$ over $\F$ \emph{(}i.e.~a Serre weight for $G(\F_p)$\emph{)} such that for all $\tau$ and $\lambda$, 
\[
Z(R_{\rbar}^{\lambda+\eta, \tau}/\varpi) = \sum_\sigma m_{\lambda,\tau}(\sigma) \cZ_\sigma(\rbar),
\]
where $m_{\lambda,\tau}(\sigma)$ denotes the multiplicity of $\sigma$ in $\ovl{\sigma(\lambda,\tau)}$.
\end{conj} 

In other words, the special fibers $R_{\rbar}^{\lambda+\eta, \tau}/\varpi$ are built out of a \emph{finite list} of basic cycles $\cZ_\sigma(\rbar)$, with multiplicities governed by the purely representation theoretic quantities $m_{\lambda,\tau}(\sigma)$. 
Conjecture \ref{conj:introBM} is known when $n=2$ and $\lambda = 0$ by work of Gee and Kisin \cite{gee-kisin}. 
When $\tau$ is a generic tame type, $m_{0,\tau}(\sigma) = 1$ for $2^f$ Serre weights $\sigma$ and is zero otherwise. 
In general, the quantities $m_{\lambda,\tau}$ are very complicated: if $\lambda=0$ and $\tau$ is tame, $m_{\lambda,\tau}(\sigma)$ computes the multiplicities of a mod $p$ Deligne--Lusztig representation, which for generic $\tau$ is given by periodic Kazhdan--Lusztig polynomials.
In particular, as the rank of $G$ grows, the special fibers $R_{\rbar}^{\lambda+\eta, \tau}/\varpi$ tend to be highly non-reduced.

As explained in \S \ref{sec:TW}, to prove Theorem \ref{thm:main}, one needs to establish Axiom \ref{ax:patch}, particularly the main bottleneck \eqref{item:fullsupport}. 
We do this  by proving Theorem \ref{thm:domain}. 
\begin{thm}\label{thm:main_def_ring}\cite[Theorem 7.3.2(2)]{MLM} Assume that $\rbar$ is \emph{semisimple} and $\tau$ is a tame inertial type which is sufficiently generic relative to $\lambda$ (in the sense of \S \ref{sec:result}). Then $R_{\rbar}^{\lambda+\eta,\tau}$ is a domain (or zero). 
\end{thm}
\begin{rmk}
\begin{enumerate}
\item
Explicit computations that suggest that Theorem \ref{thm:main_def_ring} is \emph{false} without the tameness assumption on $\rbar$ when $n > 2$ unless $n=3$ and $\lambda= 0$. 
\item If $R^{\lambda+\eta,\tau}_{\rbar} \neq 0$, then sufficient genericity of $\tau$ implies that of $\rbar$ and vice versa (generally with different choices of implicit polynomials). Because of this, the conclusion of Theorem \ref{thm:main_def_ring} also holds if we let $\rbar$ be tame and sufficiently generic but impose no genericity hypothesis on $\tau$.
\end{enumerate}
\end{rmk}
\subsection{The Emerton-Gee stack}\label{sec:EG}
 In \cite{EGstack}, Emerton--Gee constructed the moduli stack $\cX_n$ over $\Spf \cO$ of rank $n$ \'{e}tale $(\varphi,\Gamma)$-modules. By its construction, $\cX_n$ interpolates framed deformation rings in the sense that the set $\cX_n(\ovl{\F}_p)$ is in bijection with the set of continuous representations $\rbar:G_K \ra \GL_n(\overline{\F}_p)$, and framed deformation rings $R_{\rbar}$ are versal rings (in the sense of \cite[Definition 2.2.9]{EGschemetheoretic}) for $\cX_n$. Furthermore, for a Hodge--Tate cocharacter $\lambda$ and an inertial type $\tau$, they construct a $p$-flat $p$-adic formal algebraic closed substack $\cX^{\lambda,\tau}$ which is characterized by the property that its points over any finite flat $\cO$-algebra correspond to potentially crystalline representations $r$ of type $(\lambda,\tau)$. Thus $\cX^{\lambda,\tau}$ interpolates the potentially crystalline deformation rings $R_{\rbar}^{\lambda,\tau}$ as $\rbar$ varies.

The basic properties of these stacks are as follows:
\begin{thm}[Emerton--Gee]
\begin{enumerate}
\item \cite[Corollary 5.5.18]{EGstack} $\cX_n$ is a Noetherian formal algebraic stack.
\item \cite[Theorem 4.8.12]{EGstack} $\cX^{\lambda,\tau}$ is a $p$-flat $p$-adic formal algebraic stack of dimension $\dim G_{/E}/P_{\lambda/E}$.
\item \cite[Theorem 6.5.1]{EGstack} The irreducible components of the underlying reduced stack $\cX_{n,\red}$ are in bijection with the Serre weights of $G(\Fp)$. 
\end{enumerate}
\end{thm}
We let $\cC_{\sigma}$ be the irreducible component labelled by $\sigma$.
Let $\cZ_{\lambda+\eta,\tau}$ denote the top dimensional cycle of $\cX^{\lambda+\eta,\tau}/\varpi$, which has dimension independent of $\lambda$ since $\lambda+\eta$ is regular dominant. One has the following interpolation of the Breuil--M\'ezard conjecture over $\cX_n$: 
\begin{conj}[Conjecture 8.2.2 \cite{EGstack}]\label{conj:intro:BMstack}  
For each Serre weight $\sigma$, there exists an effective top-dimensional cycle $\cZ_\sigma$ on $\cX_{n,\mathrm{red}}$ such that for all $\lambda$ and inertial types $\tau$, we have
\[\cZ_{\lambda+\eta,\tau} = \sum_{\sigma} m_{\lambda,\tau}(\sigma) \cZ_\sigma.
\]
\end{conj}
\begin{rmk}\label{rmk:BM}
\begin{enumerate}
\item \label{rmk:suffmany} Conjecture \ref{conj:intro:BMstack} recovers Conjecture \ref{conj:introBM} by taking versal rings at $\rbar$. Conversely, knowledge of Conjecture \ref{conj:introBM} at \emph{sufficiently many} $\rbar$ would imply Conjecture \ref{conj:intro:BMstack}. This gives a mechanism to deduce Conjecture \ref{conj:introBM} for more general $\rbar$ from a few ``basic $\rbar$''. This allows us  to reduce Theorem \ref{thm:main}\eqref{item:genBM} to the case of semisimple $\rbar$.
\item 
 In \cite[\S 8]{MLM}, using Taylor-Wiles patching, we constructed cycles $\cZ_\sigma$ for sufficiently generic $\sigma$, which satisfies a (finite) subset of the equations postulated in Conjecture \ref{conj:intro:BMstack}. As the conjectural cycles in Conjecture \ref{conj:intro:BMstack} is expected to be compatible with Taylor-Wiles patching, the cycles constructed in \emph{loc.~cit.}~should be the ``correct'' ones.
\item (cf.~\cite[Remark 1.4.11]{MLM}) One expects $\cZ_\sigma$ to contain the irreducible component $\cC_{\sigma}$ with multiplicity one. That is, one should have a decomposition
\[\cZ_\sigma= \sum_{\sigma'} b_{\sigma',\sigma} \cC_{\sigma'}\]
with $b_{\sigma',\sigma}\geq 0$ and $b_{\sigma,\sigma}=1$. This is indeed true in the cases studied in \cite[\S 8]{MLM} and \cite{gee-kisin}. For example, in the setting of \cite{gee-kisin}, the cycles $\cZ_\sigma=\cC_{\sigma}$, unless $\sigma$ is a twist of the Steinberg weight (in particular such $\sigma$ would be non-generic), in which case $\cZ_{\sigma}$ is $\cC_{\sigma}+\cC_{\sigma'}$ for a suitable $\sigma'$ (cf.~\cite[Theorem 8.6.2]{EGstack}).
For $n>3$, it is quite difficult to compute $b_{\sigma',\sigma}$, and one does not expect $\cZ_{\sigma}=\cC_{\sigma}$ in general, even for generic $\sigma$. This is analogous to the situation of the locally analytic Breuil-M\'{e}zard conjecture studied in \cite{BHS}.
\end{enumerate}
\end{rmk}

\subsection{Local models and their geometric properties}

 Let $L\cG$ be the loop group, which is the ind-group scheme given by $L\cG(R)=\GL_n(R(\!(v+p)\!))$ for any $\cO$-algebra $R$. Consider the positive loop group scheme $L^+\cG$ over $\cO$ sending an
$\cO$-algebra 
$R$ to the subgroup of $\GL_n(R[\![v+p]\!])$ consisting of matrices that
are upper triangular mod~$v$. Note that when $p$ is invertible in $R$, $L^+\cG(R)=\GL_n(R[\![v+p]\!])$ is the positive loop group for $\GL_n$, whereas when $p=0$ in $R$, $L^+\cG(R)=\cI(R)$, the standard Iwahori group scheme.

The quotient $L^+\cG\backslash L\cG$ is represented by an ind-proper
$\cO$-ind-scheme $\Gr_{\cG}$. This is a mixed characteristic
version of the degeneration of affine Grassmannians introduced by Gaitsgory: indeed its generic fiber $\Gr_{\cG, E}$ is isomorphic to an affine Grassmannian,  
while the special fiber $\Gr_{\cG, \F}$ is isomorphic to the affine flag variety $\mathrm{Fl}$.

The affine Grassmannian has the affine Schubert stratification $\Gr_{\cG,E}=\bigcup_{\lambda}L^+ \cG_E\backslash L^+ \cG_E(v+p)^{\lambda}L^+ \cG_E$, where $\lambda$ runs over dominant coweights of $\GL_n$. Similarly, the affine flag variety $\Fl=\bigcup_{\tld{w}} \cI\backslash \cI \tld{w} \cI$, where $\tld{w}$ runs over the extended affine Weyl group $\tld{W}$.

For dominant $\lambda$, the Pappas--Zhu local model $M(\leql)$ is the Zariski closure of $L^+ \cG_E\backslash L^+ \cG_E(v+p)^{\lambda}L^+ \cG_E$ in $\Gr_{\cG}$, cf.~\cite{PZ}.   

Let $\bf{a} \in \cO^n$.  We now consider the condition 
\[
 \tag{$\star$}  v \frac{dA}{dv} A^{-1}  + A \Diag(\bf{a}) A^{-1} \in
 \left(\frac{1}{v+p}\right)\mathrm{Lie} \, L^+\cG \label{eq:monodromy}
\]
for $ A\in L\cG (R)$. This is an approximation to the monodromy condition coming from $p$-adic Hodge theory.
 This condition clearly descends to a closed condition on $\Gr_{\cG}$. 
 
\begin{defn} 
\label{defn:LMQp}
The local model $M(\lambda, \nabla_{\bf{a}})$ is
  the Zariski closure in $M(\leql)$ of the locus cut out by
    \eqref{eq:monodromy} in $L^+\cG_E\backslash L^+ \cG_E(v+p)^{\lambda}L^+ \cG_E$.
\end{defn}
Note that right multiplication by the constant diagonal torus $T^\vee$ preserves \eqref{eq:monodromy}. 
(Here, $T^\vee$ is a maximal torus in $\GL_n$ which is the dual group $G^\vee$ of the group $G = \GL_n$ which appeared in \S \ref{sec:serrewt}.)
Thus, $M(\lambda, \nabla_{\bf{a}})$ inherits a $T^\vee$-action compatible with the $T^\vee$-action on $M(\leql)$.  

By contemplating the interaction of condition \eqref{eq:monodromy} with the affine Schubert stratification, one observes:
\begin{itemize}
\item \cite[Proposition 4.1.1]{MLM} $M(\lambda,\nabla_{\bf{a}})_{/E}$ is isomorphic to $P_{\lambda}\backslash \GL_n$, hence is smooth and irreducible.
\item \cite[Theorem 4.2.4]{MLM} Provided $\bf{a}$ mod $p$ sufficiently regular, the locus cut out by \eqref{eq:monodromy} in each open Schubert cell $\cI\backslash \cI \tld{w} \cI\subset \Fl$ is an affine space, with dimension combinatorially determined by $\tld{w}$. 
\end{itemize}
Thus $M(\lambda,\nabla_{\bf{a}})$ is a degeneration of a partial flag variety, and one has control over its \emph{reduced} special fiber. 
\begin{exam}
\begin{enumerate}
\item
For example, when $n=2$ and $\lambda=(1,0)$, condition $\eqref{eq:monodromy}$ is empty, and $M(\lambda,\nabla_{\bf{a}})=M(\leql)$ is a degeneration of $\bP^1$ into a union of two $\bP^1$ crossing transversely at a point.
More generally, one has $M(\lambda,\nabla_{\bf{a}})=M(\leql)$ if and only if $\lambda$ is minuscule.
\item Suppose $n=3$ and $\lambda=(2,1,0)$, and $\bf{a}$ mod $p$ is sufficiently regular. Then $\dim M(\leql)=4$, whereas $\dim M(\lambda,\nabla_{\bf{a}})=\dim B\backslash \GL_3 =3$. The special fiber $M(\lambda,\nabla_{\bf{a}})_{\F}$ is reduced and has $9$ irreducible components, six of which are isomorphic to the flag variety $B\backslash \GL_3$, while the remaining three are more complicated rational smooth varieties. Already in this case, the behavior of the intersections among the irreducible components is somewhat elaborate, cf.~\cite{wild}. 
\end{enumerate}
\end{exam}
\begin{rmk} \label{rmk:explicit_charts}Around each point $\tld{z}\in \Fl$, one can write down an explicit open neighborhood $\cU(\tld{z})$ of $\Gr_{\cG}$ using the theory of the ``big cell''. This allows us to in principle give explicit coordinate charts for $M(\lambda,\nabla_{\bf{a}})$: the coordinate charts parametrize matrices $A$ with polynomial entries whose degrees are bounded in terms of $\tld{z}$, and one then imposes elementary divisor conditions dictated by $\lambda$ together with the explicit equation \eqref{eq:monodromy} 
and takes the $p$-saturation of the result. It is the $p$-saturation operation that makes this description rather difficult to work with.
\end{rmk}

In order to establish the connection between the above models to Galois deformation theory, we have to understand the behavior of $M(\lambda, \nabla_{\bf{a}})$ under completion. The essential difficulty is that an irreducible variety may break up into formal branches in some complicated way after completions: its singularities may not be \emph{unibranch}. 
Unfortunately, $M(\lambda, \nabla_{\bf{a}})$ fails to be unibranch in general, and in such situations it is difficult to control the subset of the formal branches that are related to Galois deformation theory. Fortunately, it turns out there is supply of special points where this difficulty does not manifest:

\begin{thm} \label{thm:intro:unibranchMLM} 
\emph{(}\cite[Theorem 3.7.1]{LLLM}\emph{)} There exists a nonzero polynomial $P \in \Z[X_1, \ldots, X_n]$ such that if $P(\bf{a}) \neq 0 \mod p$, then for any $T$-fixed point  $x \in M(\lambda, \nabla_{\bf{a}})(\ovl{\F}_p)$,  the completed local ring $\cO_{ M(\lambda, \nabla_{\bf{a}}), x}^{\wedge}$ is a domain $($i.e., $M(\lambda, \nabla_{\bf{a}})$ is unibranch at its $T$-fixed points).   
\end{thm} 
This key result, whose proof we now sketch, underlies everything else. 
One first observes that the theorem holds (under a mild assumption on the characteristic) for the equal characteristic analogues of $M(\lambda, \nabla_{\bf{a}})$ where $E$ is replaced by $\F(\!(t)\!)$. In this function field setting, there is an additional symmetry: there is an extra $\bG_m$-action given by ``loop rotation'' which scales $t$. This implies that the $T$-fixed points look like cone points, i.e.~the fixed point of an attracting torus action, and one observes that cone points are unibranch. We then deduce the mixed characteristic case by a spreading out argument. The essential point here is that unibranch can be phrased in terms of connectedness of fibers of the normalization map, and normalization is preserved by generic base change. This explains the occurence of the universal polynomial $P$: its vanishing locus is the obstruction to certain properties being preserved under base change.

\subsection{Local models and Emerton--Gee stacks}

Recall that we fixed a finite unramified extension $K/\Qp$. 
Let $k$ be the residue field of $K$. 
Let $\cJ$ be the set of embeddings $\Hom_{\Qp}(K, \overline{\Q}_p)$ which we identify with $\Hom_{\Qp}(K, E) = \Hom(k, \F)$ using the inclusion $E \subset \overline{\Q}_p$. 

To any \emph{tame} inertial type $\tau$ for $I_K$, one can associate a collection $\bf{a}_{\tau} = (\bf{a}_{\tau, j})_{j \in \cJ}$, where $\bf{a}_{\tau, j} \in \cO^n$ records the inertial weights of $\tau$. 

In the ``lowest alcove" principal series case, $\bf{a}_{\tau}$ is defined so that $\tau$ is the direct sum \[\bigoplus_{i=1}^n \prod_{j \in \cJ} j \circ \omega_K^{\bf{a}_{\tau,j}^{(i)}}\] where $\omega_K: I_K \ra k^\times$ is the reduction of the Lubin--Tate character $I_K \ra \cO_K^\times$. 
Set $\lambda = (\lambda_j)_{j \in \cJ} \in (\Z^n)^{\cJ}$, a Hodge--Tate cocharacter. Define 
\[
M_{\cJ}(\lambda, \nabla_{\bf{a}_{\tau}}) = \prod_{j \in \cJ} M(\lambda_j, \nabla_{ \bf{a}_{\tau, j}})
\] 
where, for each $j\in\cJ$, the local models $M(\lambda_j, \nabla_{ \bf{a}_{\tau, j}})$ are those appearing in  Definition \ref{defn:LMQp}.

The relationship between the local models and the Emerton--Gee stacks is given by the following:
\begin{thm}\emph{(}\cite[Theorem 7.3.2]{LLLM}\emph{)}
\label{thm:intro:LMD}   If $\tau$ is sufficiently generic \emph{(}with respect to $\lambda$\emph{)}, then
 there exist Zariski open covers $
 \underset{\text{\tiny{$\tld{z} $}}}{\text{\Large{$\bigcup$}}}
 \cX_{\reg}^{\leq \lambda,\tau}(\tld{z}) $ and $\underset{\text{\tiny{$\tld{z} $}}}{\text{\Large{$\bigcup$}}} U_{\reg}(\tld{z},\leql,\nabla_{\bf{a}_\tau})^{\wedge_p}$ of $\underset{\tiny{\substack{\lambda' \leq \lambda\\ \lambda' \text{\emph{reg. dom.}}}}}{\text{\Large{$\bigcup$}}} \cX^{\lambda',\tau} $ and  $ \underset{\tiny{\substack{\lambda' \leq \lambda\\ \lambda' \text{\emph{reg. dom.}}}}}{\text{\Large{$\bigcup$}}}M(\lambda', \nabla_{\bf{a}_\tau})^{\wedge_p}$, respectively, such that for each $\tld{z} $, there exists a local model diagram 
\begin{equation}\label{eqn:LMD}
\xymatrix{
& \tld{\cX}_{\reg}^{\leq \lambda,\tau}(\tld{z}) \ar[dl] \ar[dr] & \\
\cX_{\reg}^{\leq \lambda,\tau}(\tld{z}) & &U_{\reg}(\tld{z},\leql,\nabla_{\bf{a}_\tau})^{\wedge_p} \\
}
\end{equation}

\end{thm}

\begin{rmk} \begin{enumerate}
\item In the above diagram, the $T$-fixed points of the local models have a simple Galois theoretic interpretation: they correspond to semisimple $\rbar$.
\item 
When $\lambda = \eta$, one has
$
\underset{\tiny{\substack{\lambda' \leq \lambda\\ \lambda' \text{reg. dom.}}}}{\text{{$\bigcup$}}} \cX^{\lambda',\tau}=\cX^{\eta,\tau}.$
Since potentially crystalline deformation rings of type $(\eta,\tau)$ are versal rings to $\cX^{\eta, \tau}$, we see that they appear (up to smooth modifications) as the completion of local rings of $M(\eta, \nabla_{\bf{a}_\tau})$ at closed points. In particular, we deduce the irreducibility of the potentially crystalline deformation rings of  Theorem \ref{thm:main_def_ring} from the unibranch property of the local models (at the appropriate points). 
This completes the proof of Theorem \ref{thm:main}\eqref{item:genBM} and the first half of Theorem \ref{thm:main}\eqref{item:unitSW} (see Remark \ref{rmk:BM}(\ref{rmk:suffmany})).
\item Combining the theorem with Remark \ref{rmk:explicit_charts}, we get an algorithm to write down explicit presentations of (unions of) $R^{\lambda,\tau}_{\rbar}$ (for regular $\lambda$).
\end{enumerate}
\end{rmk}

We now give a slightly simplified outline of the proof of Theorem \ref{thm:intro:LMD}. The construction of the stacks $\cX^{\lambda,\tau}$ comes in two steps: 
\begin{itemize}
\item Using integral $p$-adic Hodge theory, one can attach Breuil--Kisin modules to lattices in (potentially) crystalline representations for $G_K$.  Thus the first step is to construct a moduli stack of Breuil--Kisin modules $Y^{\leq \lambda,\tau}$ with tame descent data of type $(\lambda,\tau)$.
\item As not all Breuil--Kisin modules come from lattices in (potentially) crystalline representations, one needs cut down $Y^{\lambda,\tau}$ by appropriate conditions to get $\cX^{\lambda,\tau}$.
\end{itemize}
Accordingly, the proof is divided into two steps:
\begin{itemize}
\item In the first step, we show that $Y^{\leq \lambda,\tau}$ is locally modelled by the Pappas-Zhu model $M(\leql)$. This is not surprising, as Breuil--Kisin modules are a projective $\cO_K[\![u]\!]$-modules with certain semi-linear structures, and thus are closely related to points of $\Gr_{\cG}$. Using the open cover $\Gr_\cG=\bigcup_{\tld{z}} \cU(\tld{z})$ (cf.~Remark \ref{rmk:explicit_charts}), we get an analogue of the local model diagram \eqref{eqn:LMD} for $Y^{\leql,\tau}$ and induced open affine covers on every object in sight. 
\item After the first step, we get \emph{two} closed substacks of $Y^{\leql,\tau}(\tld{z})$: the substack $\cX^{\leql,\tau}(\tld{z})$ and the substack $\cX^{\leql,\tau,\star}(\tld{z})$ induced by the $p$-adic completion of $ \bigcup_{\lambda' \leq \lambda} M(\lambda', \nabla_{\bf{a}_\tau})$ along the local model diagram for $Y^{\leql,\tau}$. They are genuinely different substacks, because condition (\ref{eq:monodromy}) is only the ``first order term'' of the condition cutting out $\cX^{\leql,\tau}$ inside $Y^{\leql,\tau}$. 

However, the two substacks are $p$-adically close, and using the smoothness of the generic fiber of $M(\lambda, \nabla_{\bf{a}})$, one can produce a non-canonical embedding $\cX^{\leql,\tau}(\tld{z})\into \cX^{\leql,\tau,\star}(\tld{z})$. Since both stacks turn out to have the same dimension, the maximal dimensional part $\cX_\reg^{\leql,\tau}(\tld{z})$ of $\cX^{\leql,\tau}(\tld{z})$ embeds into the maximal dimension part of $\cX^{\leql,\tau,\star}(\tld{z})$. Now, using the results of \cite{LLL} (which ultimately uses Taylor--Wiles patching, and hence automorphic forms), one obtains a lower bound on the number of irreducible components (of the spectrum of the structure sheaf) of the former, while Theorem \ref{thm:intro:unibranchMLM} gives the same upper bound for the number of irreducible components (of the spectrum of the structure sheaf) of the latter. Thus the two maximal dimension parts are (non-canonically) isomorphic to each other, which concludes the proof. 
\end{itemize}

As the above outline suggests, the arrows in the local model diagram are produced by Hensel-type lifting arguments, and thus are highly non-canonical. However, modulo $p$ this issue disappears, and the local model diagrams on the open cover produced by Theorem \ref{thm:intro:LMD} glue together. Consequently, the analysis of irreducible components of the special fibers of local models implies the following. 
\begin{thm}
\label{thm:intro:irreducible_components_mod_p} For $\tau$ sufficiently generic (with respect to $\lambda$):
\begin{enumerate}
\item $\cX^{\lambda+\eta,\tau}_\red=\cup_\sigma \cC_{\sigma}$, where the union runs over all Serre weights $\sigma\in \JH(\ovl{\sigma(\lambda,\tau)})$. 
\item There is a natural bijection between the irreducible components of $\overline{M}(\lambda+\eta, \nabla_{\bf{a}_\tau})$ and the Jordan--H\"older factors of $\ovl{\sigma(\lambda,\tau)}$.
\item For each $\sigma\in \JH(\ovl{\sigma(\lambda,\tau)})$, we have a mod $p$ local model diagram: 

\begin{equation}\label{eq:mod_p_local_model}
\xymatrix{
& \tld{\cC}_{\sigma} \ar[dl] \ar[dr] & \\
\cC_{\sigma} &   &  \overline{M}(\lambda+\eta, \nabla_{\bf{a}_\tau})_{\sigma} \\
}
\end{equation} 
where $\overline{M}(\lambda+\eta, \nabla_{\bf{a}_\tau})_{\sigma}$ is the irreducible component of $\overline{M}(\lambda+\eta, \nabla_{\bf{a}_\tau})$ labelled by $\sigma$ and both arrows are torsors for the torus $(T^\vee)^{\cJ}$ \emph{(}with respect to different $(T^\vee)^{\cJ}$-actions\emph{)}.   
\end{enumerate}
\end{thm}
\begin{rmk} 
\label{rmk:intro:MLM:sp:fib}
\begin{enumerate}
\item The proof of Theorem \ref{thm:intro:irreducible_components_mod_p} does not go through Theorem \ref{thm:intro:LMD}. Because of that it holds under much milder genericity conditions compared to our other theorems: if $\sigma(\tau) = R$, we only require that $R$ is $m$-generic for $m$ sufficiently large depending on $\lambda$ (larger than both $2\langle \lambda,\alpha^\vee\rangle + 2$ and $4n + \langle \lambda,\alpha\rangle$ for all roots $\alpha$).
\item Over $\F$, equation \eqref{eq:monodromy} becomes the equation cutting out a \emph{deformed} affine Springer fiber in $\Fl$, cf.~\cite{Frenkel_Zhu}. Thus the irreducible components $\overline{M}(\lambda+\eta, \nabla_{\bf{a}_\tau})_{\sigma}$ are irreducible components of a (product of) \emph{deformed} affine Springer fiber(s). (It is immediate from the aforementioned \cite[Theorem 4.2.4]{MLM} that these irreducible components are rational varieties.) Theorem \ref{thm:intro:irreducible_components_mod_p} then allows us to get a handle on the irreducible component $\cC_{\sigma}$ of the reduced Emerton--Gee stack for generic $\sigma$.
In particular, one can get a description of the semisimple points on $\cC_\sigma$, and this is the critical ingredient for the verification of Herzig's recipe (Theorems \ref{thm:main}\eqref{item:unitSW} and \ref{thm:main}\eqref{item:GLnSW}).
\item The irreducible components $\overline{M}(\lambda+\eta, \nabla_{\bf{a}_\tau})_{\sigma}$ are fairly easy to implement on computer algebra systems such as Macaulay2. For any given $n$, this allows us to probe the structure of $\cC_\sigma$ in a purely algorithmic manner.
\end{enumerate}
\end{rmk}

\begin{exam}
\begin{enumerate}
\item (Fontaine-Laffaille components) For $\sigma$ in the lowest alcove, the corresponding irreducible component of the deformed affine Springer fiber is isomorphic to a product of flag varieties $(B\backslash \GL_n)^{\cJ}$. We deduce from this that
\[\cC_{\sigma}=[(N\backslash \GL_n)^{\cJ}/T^{\cJ}]\]
where $N$ is the subgroup of unipotent upper triangular matrices, and $T^{\cJ}$ acts via \emph{shifted conjugation}: $(t_j)\cdot (Ng_j)=(Nt_{j}g_j t_{j\circ \phz}^{-1})$, where $\circ \phz$ denotes pre-composition with Frobenius on $K$. In particular, for $n=2$, all components $\cC_{\sigma}$ for generic $\sigma$ are of this form.
\item  
When $n=3$, there are two types of irreducible components at each factor $j\in \cJ$: the flag variety $B\backslash \GL_3$ or a more complicated rational and smooth variety (this can be extracted, for example, from the description of minimal primes in \cite[Table 3]{wild}). In particular, there are $2^f$ types of $\cC_\sigma$ (for generic $\sigma$) which correspond to the possible $p$-alcoves containing the highest weight of $\sigma$. 
\item For $n=4$, there are generic $\sigma$ for which $\cC_{\sigma}$ is singular (e.g.~for $\sigma$ with highest weight in the highest $p$-restricted alcove). Thus the smoothness of $\cC_{\sigma}$ appears to be a low rank coincidence.
\end{enumerate}
\end{exam}

\newpage
\bibliography{BiblioSurvey}

\newcommand{\etalchar}[1]{$^{#1}$}
\providecommand{\bysame}{\leavevmode\hbox to3em{\hrulefill}\thinspace}
\providecommand{\MR}{\relax\ifhmode\unskip\space\fi MR }
\providecommand{\MRhref}[2]{%
  \href{http://www.ams.org/mathscinet-getitem?mr=#1}{#2}
}
\providecommand{\href}[2]{#2}
\begin{thebibliography}{LLHLM20b}

\bibitem[ACC{\etalchar{+}}18]{10author}
Patrick~B. Allen, Frank Calegari, Ana Caraiani, Toby Gee, David Helm, Bao V.~Le
  Hung, James Newton, Peter Scholze, Richard Taylor, and Jack~A. Thorne,
  \emph{Potential automorphy over {C}{M} fields}, 2018.

\bibitem[ADP02]{ADP}
Avner Ash, Darrin Doud, and David Pollack, \emph{Galois representations with
  conjectural connections to arithmetic cohomology}, Duke Math. J. \textbf{112}
  (2002), no.~3, 521--579. \MR{1896473}

\bibitem[Bal12]{balaji}
Sundeep Balaji, \emph{G-valued potentially semi-stable deformation rings},
  ProQuest LLC, Ann Arbor, MI, 2012, Thesis (Ph.D.)--The University of Chicago.
  \MR{3152673}

\bibitem[BDJ10]{BDJ}
Kevin Buzzard, Fred Diamond, and Frazer Jarvis, \emph{On {S}erre's conjecture
  for mod {$\ell$} {G}alois representations over totally real fields}, Duke
  Math. J. \textbf{155} (2010), no.~1, 105--161. \MR{2730374 (2012k:11067)}

\bibitem[BG14]{BG}
Kevin Buzzard and Toby Gee, \emph{The conjectural connections between
  automorphic representations and {G}alois representations}, Automorphic forms
  and {G}alois representations. {V}ol. 1, London Math. Soc. Lecture Note Ser.,
  vol. 414, Cambridge Univ. Press, Cambridge, 2014, pp.~135--187. \MR{3444225}

\bibitem[BHS19]{BHS}
Christophe Breuil, Eugen Hellmann, and Benjamin Schraen, \emph{A local model
  for the trianguline variety and applications}, {Publications mathematiques de
  l' IHES} \textbf{130} (2019), 299--412.

\bibitem[BM02]{BM}
Christophe Breuil and Ariane M{\'e}zard, \emph{Multiplicit\'es modulaires et
  repr\'esentations de {${\rm GL}\sb 2({\bf Z}\sb p)$} et de {${\rm
  Gal}(\overline{\bf Q}\sb p/{\bf Q}\sb p)$} en {$l=p$}}, Duke Math. J.
  \textbf{115} (2002), no.~2, 205--310, With an appendix by Guy Henniart.
  \MR{1944572 (2004i:11052)}

\bibitem[BS73]{Borel-Serre}
A.~Borel and J.-P. Serre, \emph{Corners and arithmetic groups}, Commentarii
  mathematici Helvetici \textbf{48} (1973), 436--483.

\bibitem[Car79]{cartier}
P.~Cartier, \emph{Representations of {$p$}-adic groups: a survey}, Automorphic
  forms, representations and {$L$}-functions ({P}roc. {S}ympos. {P}ure {M}ath.,
  {O}regon {S}tate {U}niv., {C}orvallis, {O}re., 1977), {P}art 1, Proc. Sympos.
  Pure Math., XXXIII, Amer. Math. Soc., Providence, R.I., 1979, pp.~111--155.
  \MR{546593}

\bibitem[CEG{\etalchar{+}}16]{CEGGPS}
Ana Caraiani, Matthew Emerton, Toby Gee, David Geraghty, Vytautas
  Pa{\v{s}}k{\=u}nas, and Sug~Woo Shin, \emph{Patching and the {$p$}-adic local
  {L}anglands correspondence}, Camb. J. Math. \textbf{4} (2016), no.~2,
  197--287. \MR{3529394}

\bibitem[CG18]{CG}
Frank Calegari and David Geraghty, \emph{Modularity lifting beyond the
  {T}aylor-{W}iles method}, Invent. Math. \textbf{211} (2018), no.~1, 297--433.
  \MR{3742760}

\bibitem[CH13]{CH}
Ga\"{e}tan Chenevier and Michael Harris, \emph{Construction of automorphic
  {G}alois representations, {II}}, Camb. J. Math. \textbf{1} (2013), no.~1,
  53--73. \MR{3272052}

\bibitem[CS17]{CS}
Ana Caraiani and Peter Scholze, \emph{On the generic part of the cohomology of
  compact unitary {S}himura varieties}, Ann. of Math. (2) \textbf{186} (2017),
  no.~3, 649--766. \MR{3702677}

\bibitem[CV92]{CV}
Robert~F. Coleman and Jos\'{e}~Felipe Voloch, \emph{Companion forms and
  {K}odaira-{S}pencer theory}, Invent. Math. \textbf{110} (1992), no.~2,
  263--281. \MR{1185584}

\bibitem[Del71]{deligne}
Pierre Deligne, \emph{Formes modulaires et repr\'{e}sentations {$l$}-adiques},
  S\'{e}minaire {B}ourbaki. {V}ol. 1968/69: {E}xpos\'{e}s 347--363, Lecture
  Notes in Math., vol. 175, Springer, Berlin, 1971, pp.~Exp. No. 355, 139--172.
  \MR{3077124}

\bibitem[DL76]{DeligneLusztig}
P.~Deligne and G.~Lusztig, \emph{Representations of reductive groups over
  finite fields}, Ann. of Math. (2) \textbf{103} (1976), no.~1, 103--161.
  \MR{0393266}

\bibitem[Edi92]{edixhoven}
Bas Edixhoven, \emph{The weight in {S}erre's conjectures on modular forms},
  Invent. Math. \textbf{109} (1992), no.~3, 563--594. \MR{1176206}

\bibitem[EGa]{EGstack}
Matthew Emerton and Toby Gee, \emph{Moduli of \'etale
  $(\varphi,{\Gamma})$-modules and the existence of crystalline lifts}, \url{
  http://arxiv.org/abs/1908.07185}, preprint (2019).

\bibitem[EGb]{EGschemetheoretic}
\bysame, \emph{`{S}cheme theoretic images' of morphisms of stacks}, Algebraic
  Geometry, to appear.

\bibitem[EG14]{EG}
\bysame, \emph{A geometric perspective on the {B}reuil-{M}\'ezard conjecture},
  J. Inst. Math. Jussieu \textbf{13} (2014), no.~1, 183--223. \MR{3134019}

\bibitem[FM95]{fontaine-mazur}
J-M. Fontaine and B.~Mazur, \emph{Geometric {G}alois representations}, Elliptic
  curves, modular forms, \& {F}ermat's last theorem ({H}ong {K}ong, 1993), Ser.
  Number Theory, I, Int. Press, Cambridge, MA, 1995, pp.~41--78. \MR{1363495
  (96h:11049)}

\bibitem[Fra98]{franke}
Jens Franke, \emph{Harmonic analysis in weighted {$L_2$}-spaces}, Ann. Sci.
  \'{E}cole Norm. Sup. (4) \textbf{31} (1998), no.~2, 181--279. \MR{1603257}

\bibitem[FZ10]{Frenkel_Zhu}
Edward Frenkel and Xinwen Zhu, \emph{Any flat bundle on a punctured disc has an
  oper structure}, Math. Res. Lett. \textbf{17} (2010), no.~1, 27--37.
  \MR{2592725}

\bibitem[Gee11]{gee-annalen}
Toby Gee, \emph{Automorphic lifts of prescribed types}, Math. Ann. \textbf{350}
  (2011), no.~1, 107--144. \MR{2785764 (2012c:11118)}

\bibitem[GHS18]{GHS}
Toby Gee, Florian Herzig, and David Savitt, \emph{General {S}erre weight
  conjectures}, J. Eur. Math. Soc. (JEMS) \textbf{20} (2018), no.~12,
  2859--2949. \MR{3871496}

\bibitem[GK14]{gee-kisin}
Toby Gee and Mark Kisin, \emph{The {B}reuil-{M}\'ezard conjecture for
  potentially {B}arsotti-{T}ate representations}, Forum Math. Pi \textbf{2}
  (2014), e1, 56. \MR{3292675}

\bibitem[GLS14]{GLS}
Toby Gee, Tong Liu, and David Savitt, \emph{The {B}uzzard-{D}iamond-{J}arvis
  conjecture for unitary groups}, J. Amer. Math. Soc. \textbf{27} (2014),
  no.~2, 389--435. \MR{3164985}

\bibitem[GN]{GN}
Toby Gee and James Newton, \emph{Patching and the completed cohomology of
  locally symmetric spaces}, J. Inst. Math. Jussieu, to appear.

\bibitem[Gro98]{gross-satake}
Benedict~H. Gross, \emph{On the {S}atake isomorphism}, Galois representations
  in arithmetic algebraic geometry ({D}urham, 1996), London Math. Soc. Lecture
  Note Ser., vol. 254, Cambridge Univ. Press, Cambridge, 1998, pp.~223--237.
  \MR{1696481}

\bibitem[Gro99]{gross}
\bysame, \emph{Algebraic modular forms}, Israel J. Math. \textbf{113} (1999),
  61--93. \MR{1729443}

\bibitem[Her09]{herzig-duke}
Florian Herzig, \emph{The weight in a {S}erre-type conjecture for tame
  {$n$}-dimensional {G}alois representations}, Duke Math. J. \textbf{149}
  (2009), no.~1, 37--116. \MR{2541127 (2010f:11083)}

\bibitem[HT01]{HT}
Michael Harris and Richard Taylor, \emph{The geometry and cohomology of some
  simple {S}himura varieties}, Annals of Mathematics Studies, vol. 151,
  Princeton University Press, Princeton, NJ, 2001, With an appendix by Vladimir
  G. Berkovich. \MR{1876802}

\bibitem[Jan81]{jantzenDL}
Jens~Carsten Jantzen, \emph{Zur {R}eduktion modulo {$p$} der {C}haraktere von
  {D}eligne und {L}usztig}, J. Algebra \textbf{70} (1981), no.~2, 452--474.
  \MR{623819}

\bibitem[Kis08]{KisinPSS}
Mark Kisin, \emph{Potentially semi-stable deformation rings}, J. Amer. Math.
  Soc. \textbf{21} (2008), no.~2, 513--546. \MR{2373358 (2009c:11194)}

\bibitem[Kis09]{kisin-fontaine-mazur}
\bysame, \emph{The {F}ontaine-{M}azur conjecture for {${\rm GL}\sb 2$}}, J.
  Amer. Math. Soc. \textbf{22} (2009), no.~3, 641--690. \MR{2505297
  (2010j:11084)}

\bibitem[Kot92]{kottwitz}
Robert~E. Kottwitz, \emph{On the {$\lambda$}-adic representations associated to
  some simple {S}himura varieties}, Invent. Math. \textbf{108} (1992), no.~3,
  653--665. \MR{1163241}

\bibitem[KW09]{KW_Serre_2}
Chandrashekhar Khare and Jean-Pierre Wintenberger, \emph{Serre's modularity
  conjecture. {II}}, Invent. Math. \textbf{178} (2009), no.~3, 505--586.
  \MR{2551764}

\bibitem[Lab99]{labesse}
Jean-Pierre Labesse, \emph{Cohomologie, stabilisation et changement de base},
  Ast\'{e}risque (1999), no.~257, vi+161, Appendix A by Laurent Clozel and
  Labesse, and Appendix B by Lawrence Breen. \MR{1695940}

\bibitem[LHLM22]{wild}
Daniel Le, Bao Viet~Le Hung, Brandon Levin, and Stefano Morra, \emph{Serre
  weights for three-dimensional wildly ramified {G}alois representations},
  arXiv:2202.03303 (2022).

\bibitem[LLH]{LL}
Daniel Le and Bao~Viet Le~Hung, \emph{Serre weights for $\mathrm{GL}_n$ over
  {C}{M} fields}, in preparation.

\bibitem[LLHL19]{LLL}
Daniel Le, Bao~V. Le~Hung, and Brandon Levin, \emph{Weight elimination in
  {S}erre-type conjectures}, Duke Math. J. \textbf{168} (2019), no.~13,
  2433--2506. \MR{4007598}

\bibitem[LLHLM18]{LLLM}
Daniel Le, Bao~V. Le~Hung, Brandon Levin, and Stefano Morra, \emph{Potentially
  crystalline deformation rings and {S}erre weight conjectures: shapes and
  shadows}, Invent. Math. \textbf{212} (2018), no.~1, 1--107. \MR{3773788}

\bibitem[LLHLM20a]{MLM}
\bysame, \emph{Local models for {G}alois deformation rings and applications},
  arXiv:2007.05398 (2020).

\bibitem[LLHLM20b]{LLLM2}
\bysame, \emph{Serre weights and {B}reuil's lattice conjectures in dimension
  three}, Forum Math. Pi \textbf{8} (2020), e5, 135. \MR{4079756}

\bibitem[Mat67]{matsushima}
Yoz\^{o} Matsushima, \emph{A formula for the {B}etti numbers of compact locally
  symmetric {R}iemannian manifolds}, J. Differential Geometry \textbf{1}
  (1967), 99--109. \MR{222908}

\bibitem[New14]{newton13}
James Newton, \emph{Serre weights and {S}himura curves}, Proc. Lond. Math. Soc.
  (3) \textbf{108} (2014), no.~6, 1471--1500. \MR{3218316}

\bibitem[PZ13]{PZ}
George Pappas and Xinwen Zhu, \emph{Local models of {S}himura varieties and a
  conjecture of {K}ottwitz}, Invent. Math. \textbf{194} (2013), no.~1,
  147--254. \MR{3103258}

\bibitem[Sch15]{scholze}
Peter Scholze, \emph{On torsion in the cohomology of locally symmetric
  varieties}, Ann. of Math. (2) \textbf{182} (2015), no.~3, 945--1066.
  \MR{3418533}

\bibitem[Ser77]{serre-book}
Jean-Pierre Serre, \emph{Linear representations of finite groups},
  Springer-Verlag, New York-Heidelberg, 1977, Translated from the second French
  edition by Leonard L. Scott, Graduate Texts in Mathematics, Vol. 42.
  \MR{0450380}

\bibitem[Ser87]{serre-duke}
\bysame, \emph{Sur les repr\'esentations modulaires de degr\'e {$2$} de {${\rm
  Gal}(\overline{\bf Q}/{\bf Q})$}}, Duke Math. J. \textbf{54} (1987), no.~1,
  179--230. \MR{885783 (88g:11022)}

\bibitem[Shi11]{shin}
Sug~Woo Shin, \emph{Galois representations arising from some compact {S}himura
  varieties}, Ann. of Math. (2) \textbf{173} (2011), no.~3, 1645--1741.
  \MR{2800722}

\bibitem[Zhu14]{Zhu_coherence}
Xinwen Zhu, \emph{On the coherence conjecture of {P}appas and {R}apoport}, Ann.
  of Math. (2) \textbf{180} (2014), no.~1, 1--85. \MR{3194811}

\end{thebibliography}
\bibliographystyle{amsalpha}

\end{document}